\documentstyle[11pt,amssymb,amsfonts]{article}

\begin{document}
\newcommand{\p}{\parallel }
\makeatletter \makeatother
\newtheorem{th}{Theorem}[section]
\newtheorem{lem}{Lemma}[section]
\newtheorem{de}{Definition}[section]
\newtheorem{rem}{Remark}[section]
\newtheorem{cor}{Corollary}[section]
\renewcommand{\theequation}{\thesection.\arabic {equation}}

\title{{\bf Volterra calculus, local equivariant family index theorem and
equivariant eta forms}
}
\author{ Yong Wang }

\date{}
\maketitle

\begin{abstract} In this paper, we give proofs of the family index formula and the equivariant family index formula
by the Greiner's approach to heat kernel asymptotics. We compute
equivariant family JLO characters.
We also define the equivariant eta form and give a proof of its regularity.\\

\end{abstract}

\section{Introduction}
    \quad The first success of proving the Atiyah-Singer index theorem directly by heat kernel method
was achieved by Patodi [Pa], who carried out the ``fantastic
cancellation'' (cf. [MS]) for the Laplace operators and who for the
first time proved a local version of the Gauss-Bonnet-Chern theorem.
After that
 several different direct heat kernel proofs of the Atiyah-Singer index theorem for
Dirac operators appeared independently: Bismut [Bi1], Getzler
[Ge1], [Ge2] and Yu [Yu], Ponge [Po1]. All the proofs have their own
advantages. Motivated by the problem of generalizing the heat kernel
proofs of the index theorem to prove a local index theorem for
families of elliptic operators, Quillen [Qu] introduced the concept
of superconnections, which was developed by Bismut to give a heat
kernel representation for the Chern character of families of first
order elliptic operators. Then using his probabilitistic method,
Bismut [Bi2] obtained a proof of the local index theorem for families
of Dirac operators. In [BV1],[Do1], they gave two different proofs of
local index theorem for families of Dirac operators. Using the
method of Yu, Zhang gave another proof of the local index theorem
for families of Dirac operators in [Zh1]. The first purpose of this paper is
to give another proof of
   the local index
theorem for families of Dirac operators by the Ponge's method in [Po1].\\
\indent The Atiyah-Bott-Segal-Singer index formula is a
generalization of the Atiyah-Singer index theorem to the case with group actions.
In [BV2], Berline and Vergne gave a heat kernel proof of the
Atiyah-Bott-Segal-Singer index formula. In [LYZ], Lafferty, Yu and
Zhang presented a very simple and direct geometric proof for
equivariant index of the Dirac operator. In [PW], Ponge and Wang
gave a different proof of the equivariant index formula by the
Greiner's approach of the heat kernel asymptotics. In [LM], in order
to prove family rigidity theorems, Liu and Ma proved the equivariant
family index formula.
 The second  purpose of this paper is to give another proof of
   the local equivarint index
theorem for families of Dirac operators by the Greiner's approach of the heat kernel asymptotics.  \\
     \indent It is known, due to Connes [Co] and an
equivalent but convenient version due to Jaffe, Lesniewski
and Osterwalder, also known as a JLO formula [JLO], that the Chern
character of a $\theta$-summable Fredholm module $(H,D)$ over a unital
$C^*$-algebra $A$, takes value in the entire cyclic cohomology of $A$.
JLO characters were computed in [CM1] and [BlF].
An explicit formula of an equivariant JLO character, associated
to the invariant Dirac operator, in the presence of a countable
discrete group action on a smooth compact spin Riemannian manifold,
was given by Azmi, [Az1]. Moreover, in [Az1] it was shown that
this equivariant cocycle is an element of the delocalized cohomology,
and it can be paired with an equivariant $K$-theory idempotent. In the case when $G$
is a compact Lie group, Chern and Hu [CH] gave an explicit formula
of an equivariant Chern-Connes character, associated to a G-equivariant
$\theta$-summable Fredholm module. In [PW], Ponge and Wang computed equivariant Chern-Connes characters
by the Greiner's approach of the heat kernel asymptotics.\\
\indent
On the other hand, Wu [Wu] constructed a bivariant Chern-Connes character for (a special class of)
$\theta$ summable modules, by incorporating the JLO formula and the superconnectiom
formalism of Quillen. Wu's bivariant character takes values in the bivariant cyclic theory described by
Lott [Lo], who constructed it as a combination of entire cyclic (co)homology and
noncommutative de Rham homology of graded differential algebra. Then
by adopting Wu's method and employing Bismut's superconnection together with
the canonical order calculus developed by Simon in [CFKB], Azmi [Az2] expressed the local formula for families in terms of
differential forms on the base and the Chern roots of the fibration.
In [Zh2], Zhang announced a different method to compute the family JLO characters by developing the family version of the method in [CH] and [Fe].
In the third part of this paper, we give details of the announcement in [Zh2] and compute the equivariant family JLO characters.\\
 \indent
    In [APS], Atiyah, Patodi and Singer proved the Atiyah-Patodi-Singer index theorem
    for manifolds with boundary and they introduced the eta
    invariants. Bismut and Freed gave a simple proof of the regularity of
    eta invariants in [BiF]. In [Po2], Ponge gave another proof of the regularity of
    eta invariants using the method in [Po1].
     Bismut and Cheeger generalized the Atiyah-Patodi-Singer index
    theorem to the family case in [BC1,2]. They used the eta form
    for families of Dirac operators. The regularity of the eta form was
    proved by the probabilistic method in [BiGS]. Donnelly generalized the Atiyah-Patodi-Singer index theorem
    to the equivariant case and introduced the equivariant eta
    invariant in [Do2]. Zhang proved the regularity of the equivariant eta
    invariant by the Clifford asymptotics in [Zh3]. In this paper, we
    firstly prove the regularity of the equivariant eta invariant by a similar method in [PW] and by introducing the Grassmann variable.
     Then we define the equivariant eta form and prove
   its regularity.\\
   \indent Using the approach of Ponge and Ponge-Wang to give new proofs of the equivariant family index theorem and the regularity of equivariant eta forms has
   two advantages. One is using the Volterra pseudodifferential calculus to get the heat kernel asymtotic expansion instead of
   heat equation discussions as in [Yu], [Zh1], [BGV, Chapter 2]. In [LYZ] and [Zh3], for proving the equivariant local index theorem and the
   regularity of the equivariant eta invariants, transformed formulas between normal and tubular coordinates are needed. In [LM], Liu and Ma used the finite propagation speed method in [Bi3] to prove the equivariant family index theorem. The other advantage is that the transformed formulas between normal and tubular coordinates are the consequence of a standard change of variable formula for pseudodifferential symbols.\\
    \indent
    This paper is organized as follows: In Section 2, we give another proof of
   the local index
theorem for families of Dirac operators by the Ponge's method. In Section 3, we give another proof of
   the local equivarint index
theorem for families of Dirac operators by the Greiner's approach of the heat kernel asymptotics. In Section 4,
we compute the equivariant family JLO characters. In Section 5, we define the equivariant eta form and prove
    its regularity.

\section{The local index
theorem for families of Dirac operators
 }
{\bf 2.1} {\bf The Greiner's approach of heat kernel asymptotics with form coefficients.}\\

 \indent In [Gr],[BGS], Greiner and Beals-Greiner-Stanton defined Volterra pseudodifferential operators and gave the heat kernel asymptotics by
 Volterra pseudodifferential calculus. A good summary on the Greiner's approach of the heat kernel asymptotics was given in [Po1]. In [Po1], Ponge
 gave a proof of the local index theorem by the Greiner's approach of the heat kernel asymptotics. In this section, we shall extend Volterra pseudodifferential calculus to the case with form coeffcients, i.e. the family case. A local family index theorem was proved originally by Bismut in
[Bi2]. A superconnection due to Bismut played an important role in the proof of the local family index theorem.\\
 \indent Let us recall the definition of superconnection due to Bismut.
 Let $M$ be an $n+\overline{q}$ dimensional compact connected manifold
 and $B$ be a $\overline{q}$ dimensional compact connected
 manifold. We assume that $\pi :M\rightarrow B$ is a submersion of
 $M$ onto $B$, which defines a fibration of $M$ with fibre $Z$. For
 $y\in B$, $\pi^{-1}(y)$ is then a submanifolds $M_y$ of $M$. Denote by $TZ$
 the $n$-dimensional vector bundle on $M$ whose fibre $T_xM_{\pi x}$
 is the tangent space at $x$ to the fibre $M_{\pi (x)}$. We assume
 that $M$ and $B$ are oriented. We take a smooth
 horizontal subbundle $T^HM$ of $TM$. Vector fields
 $X\in \Gamma (B,TB)$ will be identified with their horizontal lifts $X\in  \Gamma (M, T^HM)$. Moreover $T^H_xM$ is isomorphic to $T_{\pi(x)}B$ via
 $\pi_*$. We take a Riemannian metric on $B$ and then lift the
 Euclidean scalar product $g_B$ of $TB$ to $T^HM$.
We further assume that $TZ$ is endowed with a scalar product $g_Z$. Thus
we can introduce on $TM$ a new scalar product $g_B\oplus g_Z$, and
denote by $\nabla^L$ the Levi-Civita connection on $TM$ with respect
to this metric. Set $\nabla^B$ denote the Levi-Civita connection on
$TB$ and we still denote by $\nabla^B$ the pullback connection on
$T^HM$. Let $\nabla^Z=P_Z(\nabla^L)$ where $P_Z$ denotes the orthogonal
projection to $TZ$. Set $\nabla^{\oplus}=\nabla^B\oplus \nabla^Z$
and $S=\nabla^L-\nabla^{\oplus}$ and $T$ be the torsion tensor of
$\nabla^{\oplus}$. Denote by $SO(TZ)$ the $SO(n)$ bundle of oriented
orthonormal frames in $TZ$. Now we assume that bundle $TZ$ is spin.
Denote by $S(TZ)$ the associated spinor bundle and $\nabla^Z$ can be
lifted to a connection on $S(TZ)$. Let $D$ be the Dirac
operator in the tangent direction defined by $D=\sum_{i=1}^nc(e^*_i)\nabla^{S(TZ)}_{e_i}$ where $\nabla^{S(TZ)}$ is a spin connection on $S(TZ)$. Set $E$ be the vector bundle $\pi^*(\wedge T^*B)\otimes S(TZ)$. This bundle
carries a natural action $m_0$ of the degenerate Clifford module denoted as $Cl_0(M)$. The Clifford action of a
horizontal cotangent vector
$\alpha\in\Gamma(M,T_H^*M)$ is given by the exterior multiplication $m_0(\alpha)=\varepsilon(\alpha)$
acting on the first factor $\bigwedge T_H^*M$ in $E$, while the Clifford action of a vertical cotangent vector simply
is its Clifford action on $S(TZ)$. Define the connection by ([BGV])
$$\nabla^{E,\oplus}:=\pi^*\nabla^B\otimes 1+1\otimes \nabla^{S},\eqno(2.1)$$
$$\omega(X)(Y,Z):=g(\nabla^L_XY,Z)-g(\nabla^{\oplus}_XY,Z),\eqno(2.2)$$
$$\nabla^{E,0}_X:=\nabla^{E,\oplus}_X+\frac{1}{2}m_0(\omega(X)),\eqno(2.3)$$
for $X,Y,Z\in \Gamma (M,TM)$. Then the Bismut connection acting on $\Gamma(M,\wedge(\pi^*T^*B)\otimes S(TZ))$ is defined by
$${\mathcal{B}}=\sum_{i=1}^nc(e_i^*)\nabla_{e_i}^{E,0}+\sum_{j=1}^{\overline{q}}c(f_j^*)\nabla_{f_j}^{E,0}.\eqno(2.4)$$
Then by Proposition 10.15 in [BGV],
${\mathcal{B}}=D+\textrm{A}_{[+]},$ where $\textrm{A}_{[+]}$ is an
operator with coefficients in $\Omega_{\geq 1}(B)$. By Theorem 10.17
in [BGV], we have
$$ \textrm{F}={\mathcal{B}}^2=-\sum_{i=1}^n(\nabla_{e_i}^{E,0})^2+\sum_{i=1}^n\nabla^{E,0}_{\nabla^{TZ}_{e_i}{e_i}}+\frac{1}{4}r
=D^2+F_{[+]},\eqno(2.5)$$ where $r$ is the scalar curvature of
fibres. Let $\triangle_k$ be the $k$-simplex defined by
$$\{(\sigma_0,\cdots,\sigma_k)|\sigma_0+\cdots \sigma_k=1,~0\leq \sigma_j\leq 1\}.$$
 For a fixed $t>0$, define the operator $e^{-t{F}}$ by
$$e^{-t{F}}=e^{-tD^2}+\sum_{k>0}(-t)^kI_k,\eqno(2.6)$$
where
$$I_k=\int_{\triangle_k}e^{-\sigma_0tD^2}F_{[+]}e^{-\sigma_1tD^2}F_{[+]}\cdots e^{-\sigma_{k-1}tD^2}F_{[+]}e^{-\sigma_ktD^2}d\sigma,\eqno(2.7)$$
and where the sum is finite (see [BGV, p. 312]). Then
$$(\frac{\partial}{\partial t}+F)e^{-t{F}}=0,~~Fe^{-t{F}}=e^{-t{F}}F.\eqno(2.8)$$
In the following, we formulate the Greiner's approach of heat kernel
asymptotics with coefficients in $\wedge T_z^*B$ as in [Gr] and
[BeGS]. Let $Q_0$ be given by
 $$(Q_0u)(x,s)=\int_{0}^{\infty}e^{-sF}[u(x,t-s)]dt,~~u\in \Gamma_c(M_z\times {\mathbb{R}},S(TM_z))\otimes \wedge T_z^*B.\eqno(2.9)$$
The operator $Q_0$ maps continuously from $u$ to $D'(M_z\times {\mathbb{R}},S(TM_z))\otimes \wedge T_z^*B$ which is the dual space of $\Gamma_c(M_z\times {\mathbb{R}},S(TM_z))\otimes \wedge T_z^*B.$ By (2.8), we have
 $$(F+\frac{\partial}{\partial t})Q_0u=Q_0(F+\frac{\partial}{\partial t})u=u,~~~u\in \Gamma_c(M_z\times {\mathbb{R}},S(TM_z))\otimes \wedge T_z^*B,\eqno(2.10)$$
We define the operator
$$Q=(F+\frac{\partial}{\partial t})^{-1}=(D^2+\frac{\partial}{\partial t})^{-1}+\sum_{k>0}(-1)^k(D^2+\frac{\partial}{\partial t})^{-1}
[F_{[+]}(D^2+\frac{\partial}{\partial t})^{-1}]^k,\eqno(2.11)$$
where $(D^2+\frac{\partial}{\partial t})^{-1}$ is the Volterra inverse of $D^2+\frac{\partial}{\partial t}$ as in [BeGS]. Then
$$(F+\frac{\partial}{\partial t})Q=I-R_1;~~Q(F+\frac{\partial}{\partial t})=1-R_2,\eqno(2.22)$$
where $R_1,R_2$ are smoothing operators. Let $K_{Q_0}(x,y,t-s)$ be the distribution kernel of $Q_0$. That is
$$(Q_0u)(x,t)=\int_{M_z\times {\mathbb{R}}}K_{Q_0}(x,y,t-s)u(y,s)dyds,\eqno(2.23)$$
and let $k_t(x,y)$ be the heat kernel of $e^{-tF}$. Similar to the
discussions in [BeGS p.363], we get
$$K_{Q_0}(x,y,t)=\left\{\begin{array}{lcr}k_t(x,y)~& {\rm when}~ t>0,
\\
~ 0 & {\rm when }~ t<0.
\end{array}\right.
\eqno(2.24)$$\\

\noindent{\bf Definition 2.1}~~{\it The operator $P$ is called a Volterra $\Psi DO$ if
 (i) $P$ has the Volterra property, i.e. it has a distribution
kernel of the form $K_P(x,y, t-s)$ where $K_P(x, y, t)$ vanishes on the region $t < 0.$}\\
(ii) {\it The heat operator $P +
\frac{\partial}{\partial t}$ is parabolically homogeneous, i.e. the homogeneity with respect to
the dilations of ${\mathbb{R}}^n\times{\mathbb{R}} ^1$ is given by}
 $$\lambda\cdot (\xi,\tau)=(\lambda\xi,\lambda^2\tau),~~~~~~~(\xi,\tau)\in {\mathbb{R}}^n\times {\mathbb{R}}^1,~~\lambda\neq 0.\eqno(2.25)$$\\

\indent In the sequel for $g\in  {\mathcal{S}} ({\mathbb{R}}^{n+1})$
and $\lambda\neq 0$, we let $g_{\lambda}$ be the tempered
distribution defined by
$$\left<g_\lambda(\xi,\tau),u(\xi,\tau)\right>=|\lambda|^{-(n+2)}\left<g(\xi,\tau),u(\lambda^{-1}\xi,\lambda^{-2}\tau)\right>,~~u\in
{\mathcal{S}} ({\mathbb{R}}^{n+1}).\eqno(2.26)$$\\

\noindent{\bf Definition 2.2}~~A distribution $ g\in {\mathcal{S}}
({\mathbb{R}}^{n+1})$ is parabolic homogeneous of degree $m$, where $m \in {\mathbb{Z}}
,$ if for
any $\lambda\neq 0$, we have $g_\lambda = \lambda^m g.$\\

\indent Let ${\mathbb{C}}_-$ denote the complex halfplane $\{{\rm Im}\tau < 0\}$ with closure $\overline{{\mathbb{C}}_-}$. Then:\\

\noindent {\bf {Lemma 2.3}} ([BeGS, Prop. 1.9]). {\it Let $q(\xi, \tau)\in C^{\infty}(({\mathbb{R}}^n\times {\mathbb{R}})\backslash 0)\otimes \wedge T_z^*B$ be a parabolic homogeneous symbol
of degree $m$ such that:}\\
(i){\it ~ $q$ extends to a continuous function on $(
{\mathbb{R}}^n\times \overline{{\mathbb{C}}_-})\backslash 0 $ in such way to
be holomorphic in the
last variable when the latter is restricted to ${{\mathbb{C}}}_-$.}\\
{\it Then there is a unique $g\in  {\mathcal{S}} ({\mathbb{R}}^{n+1})\otimes \wedge T_z^*B$ agreeing with q on ${\mathbb{R}}^{n+1}\backslash 0$ so that:}\\
(ii) {\it $g$ is homogeneous of degree $m$};\\
(iii) {\it The inverse Fourier transform $\check{g}(x, t)$ vanishes for $t < 0.$}\\

\indent Let $U$ be an open subset of ${\mathbb{R}}^n$. We define Volterra symbols and Volterra $\Psi DO$¡¯s on $U\times {\mathbb{R}}^{n+1}/0$
as follows.\\

\noindent {\bf {Definition 2.4}}~~{\it The set $S_V^m(U\times
{\mathbb{R}}^{n+1})\otimes \wedge T^*_zB,~m\in{\mathbb{Z}}$ ,
consists of smooth functions $q(x, \xi, \tau)$ on $U\times
{\mathbb{R}}^{n}\times {\mathbb{R}}$ with
an asymptotic expansion $q\sim \sum_{j\geq 0}q_{m-j},$ where:}\\
-{\it  $q_l\in C^{\infty}(U\times [({\mathbb{R}}^{n}\times {\mathbb{R}})\backslash 0])\otimes \wedge T^*_zB$ is a homogeneous Volterra symbol of degree $l$, i.e. $q_l$ is parabolic
homogeneous of degree $l$ and satisfies the property (i) in Lemma 2.3 with respect to the last $n + 1$
variables;}\\
- {\it The sign $\sim$ means that, for any integer $N$ and any compact $K\subset U,$ there is a constant
$C_{NK\alpha\beta k}>0$ such that for $x\in K$ and for $|\xi|+|\tau|^{\frac{1}{2}}>1$ we have}

$$||\partial_x^\alpha\partial_\xi^\beta\partial_\tau^k(q-\sum_{j<N}q_{m-j})(x,\xi,\tau)||\leq
C_{NK\alpha\beta k}(|\xi|+|\tau|^{\frac{1}{2}})^{m-N-|\beta|-2k}.\eqno(2.27)$$
{\it For $q=\sum_lq_l\omega^l$ where $q_l\in S_V^m(U\times
{\mathbb{R}}^{n+1})$ and $\omega^l\in \wedge^l T^*_zB$, we define
$||q||=\sum_l|q_l|||\omega^l||$ and $||\omega^l||$ is the norm of $\omega^l$ in  $(\wedge^l T^*_zB,g^{TB}_z).$}\\

\noindent {\bf Definition 2.5}  {\it The set $\Psi_V^m(U\times
{\mathbb{R}},\wedge T^*_zB), ~m\in{\mathbb{Z}}$ , consists of
continuous operators $Q$ from $C_c^{\infty}(U_x\times
{\mathbb{R}}_t,\wedge T^*_zB)$ to $C^{\infty}(U_x\times
{\mathbb{R}}_t,\wedge T^*_zB)$
such that:}\\
(i) {\it $ Q $ has the Volterra property;}\\
(ii) {\it $Q = q(x,D_x,D_t) + R$ for some symbol $q$ in $S_V^m(U\times {\mathbb{R}},\wedge T^*_zB)$
and some smoothing operator $R$.}\\

\indent In the sequel if $Q$ is a Volterra $\Psi DO$, we let $K_Q(x, y, t-s)$ denote its distribution kernel, so that
the distribution $K_Q(x, y, t)$ vanishes for $t< 0$.\\

\noindent {{\bf Definition 2.6} {\it  Let $q_m(x, \xi, \tau)\in
C^{\infty}(U\times ({\mathbb{R}}^{n+1}\backslash 0))\otimes\wedge T^*_zB$ be a
homogeneous Volterra symbol of order $m$ and let $g_m \in
C^{\infty}(U)\otimes {\mathcal{S}}'({\mathbb{R}}^{n+1})\otimes\wedge
T^*_zB$ denote its unique homogeneous extension given by Lemma 2.3.
Then:}\\
- {\it $\breve{q}_m(x, y, t)$ is the inverse Fourier transform of $g_m(x, \xi, \tau )$ in the last $n + 1$ variables;}\\
- {\it $q_m(x,D_x,D_t)$ is the operator with kernel $\breve{q}_m(x, y-x, t).$}\\

The composition of ordinary Volterra symbols naturally extends to a composition
of Volterra symbols with form coefficients. Let ${\bf Q}$ and ${\bf Q'}$ be in $\Psi_V^{m_1}(U\times
{\mathbb{R}},\wedge^j T^*_zB)$ and $\Psi_V^{m_2}(U\times
{\mathbb{R}},\wedge^l T^*_zB)$ with the symbols ${\bf q}$ and ${\bf q'}$. Let the composition of ${\bf Q}$ and ${\bf Q'}$ have
the symbol ${\bf q}\widetilde{\circ} {\bf q'}$, then
$${\bf q}\widetilde{\circ} {\bf q'}=\omega_1 \wedge \omega_2\otimes q\circ q'\eqno(2.28)$$
where ${\bf q}=\omega_1 q ,~{\bf q'}=\omega_2 q'$ and
$q\circ q'$ is the ordinary composition of symbols corresponding to the
   Volterra $\Psi DO$ algebra multiplication ([Gr, BeGs]). Thus we have
$$\Psi_V^{m_1}(U\times
{\mathbb{R}},\wedge^j T^*_zB)\times \Psi_V^{m_2}(U\times
{\mathbb{R}},\wedge^l T^*_zB)\rightarrow \Psi_V^{m_1+m_2}(U\times
{\mathbb{R}},\wedge^{j+l} T^*_zB).\eqno(2.29)$$\\

\noindent{\bf Proposition 2.7}~~{\it The following properties hold.\\
1) Composition. Let $Q_j\in \Psi_V^{m_j}(U\times {\mathbb{R}})\otimes\wedge T^*_zB, ~j=1,2$
have symbol $q_j$ and suppose that $Q_1$ or $Q_2$ is
properly supported (see {\rm [Ta, p. 43]}). Then $Q_1Q_2$ is a Volterra $\Psi DO$ of order $m_1+m_2$ with symbol $q_1\widetilde{\circ }q_2\sim \sum\frac{1}{\alpha!}\partial^\alpha_\xi
q_1D^\alpha_xq_2.$}\\
2) {\it Parametrices. Let $Q=Q^m+Q^{<m}$ where $Q^m$ is in $\Psi_V^{m}(U\times {\mathbb{R}})$(an
order $m$ Volterra $\Psi DO$ without form coefficients in $\wedge T_z^*B$) and $Q^{<m}$ is in $\Psi_V^{<m}(U\times {\mathbb{R}})\otimes\wedge T^*_zB$. We assume that there is an operator $P$ in $\Psi_V^{-m}(U\times {\mathbb{R}})$ such that
$$Q^mP=1-R_1,~~~PQ^m=1-R_2\eqno(2.30)$$
where $R_1,~R_2$ are smoothing operators in $\Psi_V^{-\infty}(U\times {\mathbb{R}}).$
 let $$\widetilde{Q}=P+\sum_{k>0} (-1)^kP[Q^{<m}P]^k,\eqno(2.31)$$
then
$$Q\widetilde{Q}=1-\widetilde{R_1},~~~\widetilde{Q}Q=1-\widetilde{R_2}\eqno(2.32)$$
where $\widetilde{R}_1,~\widetilde{R}_2$ are smoothing operators in $\Psi_V^{-\infty}(U\times {\mathbb{R}})\otimes\wedge T^*_zB$.}\\

\noindent{\bf Proof.} The claim 1) comes from (2.28) and 1) of Proposition 1 in [Po1]. By (2.30) and (2.31) and direct computations, we get (2.32).
~~~~~~~~~~~~$\Box$\\

\indent By (2.11) and the fact that $(D^2+ \partial_t)^{-1}$ is a Volterra $\Psi DO$ of order $-2$ and that $F_{[+]}$ is a first order Volterra $\Psi DO$, we get\\

\noindent{\bf Proposition 2.8}~{\it The differential operator $F +
\partial_t$ is invertible and
its inverse $(F + \partial_t)^{-1}$ is a Volterra $\Psi DO$ of order $-2$.}\\

Let $Q\in \Psi_V^m(U\times {\mathbb{R}},
\wedge T^*_zB)$ have symbol $q\sim\sum q_{m-j}$ and
$q_{m-j}=\sum_{s=1}^{2^{\overline{q}}}q_{m-j,s}\omega^s$ where $q_{m-j,s}\in S^{m-j}(U\times {\mathbb{R}})$ and $\omega^s\in \wedge^sT^*B.$ We define
the inverse Fourier transform of $q_{m-j}$ by  $\check{q}_{m-j}=\sum_{s=1}^{2^{\overline{q}}}\check{q}_{m-j,s}\omega^s.$ Then we have\\

\noindent {\bf Lemma 2.9} (Compare with Lemma 2 in [Po1])~ Let $Q\in \Psi_V^m(U\times {\mathbb{R}},
\wedge T^*_zB)$, we have in local coordinates
$$K_Q(x,x,t)\sim t^{-(\frac{n}{2}+[\frac{m}{2}]+1)}\sum_{l\geq 0}t^l \check{q}_{2[\frac{m}{2}]-2l}(x,0,1),\eqno(2.33)$$
\\
\noindent{\bf Proof.} Let $Q=\sum_{r=1}^{2^{\overline{q}}}Q_r\omega_r$ where $Q_r\in \Psi_V^m(U\times {\mathbb{R}})$ and $\omega_r\in  \wedge^rT^*B.$
We note that the leading symbol of $Q_r$ is probably zero. We set that the symbol $q^{[r]}$ of $Q_r$ has an asymptotic expansion $q^{[r]}\sim
\sum_{j\geq 0} q^{[r]}_{m-j}.$
 By Lemma 2 in [Po1], we have
$$K_{Q_r}(x,x,t)\sim t^{-(\frac{n}{2}+[\frac{m}{2}]+1)}\sum_{l\geq 0}t^l \check{q}^{[r]}_{2[\frac{m}{2}]-2l}(x,0,1),\eqno(2.34)$$
By the equality
$$K_{Q}(x,x,t)=\sum_{r=1}^{2^{\overline{q}}}K_{Q_r}(x,x,t)\omega_r,\eqno(2.35)$$
 and the definition of the inverse Fourier transform of $q_{m-j}$, we get (2.33).~~~~~~~$\Box$\\

By (2.24) and Proposition 2.8 and Lemma 2.9, we get\\

\noindent {\bf Theorem 2.10} (Compare with Theorem 1.6.1 in [Gr]) {\it In
$C^{\infty}(M_z, {\rm End }(S(T(M_z))))\otimes \wedge T^*_zB$, we
have}
$$k_t(x,x)\sim t^{-\frac{n}{2}}\sum_{l\geq 0}t^l a_l(F)(x), ~t\rightarrow 0^+,~~a_l(F)(x)=\check {q}_{-2-2l}(x,0,1),\eqno(2.36)$$
{\it where ${q}_{-2-2l}(x,\xi,\tau)$ is the $-2l-2$ order symbol of $(\partial_t+F)^{-1}$ and the second equality in (2.36) holds in local
coordinates.}\\

By the same reason with Lemma 2.9 and Proposition 2 in [Po1], we have\\

\noindent {\bf Proposition 2.11} (Compare with Proposition 2 in [Po1]) {\it Let $P :
C^{\infty}(M_z, S(TM_z))\rightarrow C^{\infty}(M_z, S(TM_z))$ be a
differential operator of order $m$ and let $h_t(x, y)$ denote the
distribution kernel of $Pe^{-tF}$. Then in $C^{\infty}(M_z, {\rm End
}(S(TM_z)))\otimes \wedge T^*_zB$, we have}
$$h_t(x,x)\sim t^{-([\frac{m}{2}]+\frac{n}{2})}\sum_{l\geq 0}t^l b_l(F)(x), ~t\rightarrow 0^+,~~b_l(x)=\check{q}_{2[\frac{m}{2}]-2-2l}(x,0,1),\eqno(2.37)$$
{\it where ${q}_{2[\frac{m}{2}]-2-2l}(x,\xi,\tau)$ is the
$2[\frac{m}{2}]-2l-2$ order symbol of $P(\partial_t+F)^{-1}$ and the second equality in (2.37) holds in local
coordinates.}\\

\noindent {\bf 2.2 The local family index formula.}\\

\indent In [Bi2],[BV1],[Do1],[Zh1], several different
proofs of the local index theorem for families of Dirac operators were given. In this section, we shall give a new proof of the local family
index formula by using the Greiner's approach of the heat kernel
asymptotics with form coefficients. Comparing with previous proofs, we use the Volterra calculus with form coefficients instead of heat equation
 discussions to get the family heat kernel asymototics as in [Po1].\\
\indent Let us introduce some notations. For $z\in B$, denote by $D_z$
the restriction of D to the fiber acting on $\Gamma(M_z,S(TM_z))$. We set that the dimension $n$ of fibre to be even. Then $\Gamma(M_z,S(TM_z))$ has a splitting as a sum of $\Gamma^+(M_z,S(TM_z))$ and
$\Gamma^-(M_z,S(TM_z)).$ The operator $
D_z$ interchanges $\Gamma^+(M_z,S(TM_z))$ and
$\Gamma^-(M_z,S(TM_z))$. Let $D_{z,+},~D_{z,-}$ be the restrictions of
$D_z$ to $\Gamma^+(M_z,S(TM_z))$ and $\Gamma^-(M_z,S(TM_z))$ respectively. By Chapter 9 in [BGV], the difference bundle $[{\rm ker}D_{z,+}]
-[{\rm ker}D_{z,-}]$ over $B$ is well defined in the sense of $K$-theory. The
family index theorem presents a calculation of the Chern
character of the difference bundle as a differential form over $B$
explicitly. We change the normalization constant in the definition
of the Chern character. Namely, for a vector bundle $V$ with
connection form $\gamma$ and curvature $ C$, we set ${\rm Ch}(V) =
{\rm
Tr}({\rm exp}(-C))$. Let $\widehat{A}(R^{TZ})={\rm det}^{\frac{1}{2}}\left(\frac{R^{TZ}/2}{{\rm sinh}(R^{TZ}/2)}\right)$
and $\int_{M_z}$ denote the integral along the fibre.
 Then we have\\

\noindent{\bf Theorem 2.12} (Atiyah-Singer [AS]) {\it It holds that the form}
$$(2i\pi)^{-\frac{n}{2}}\int_{M_z}\widehat{A}(R^{TZ})\eqno(2.38)$$
{\it is a representative of ${\rm {Ch}}([{\rm ker}D_{z,+}]-[{\rm
ker}D_{z,-}])$.}\\

Let $\overline{F}$ be a complex vector bundle on $M$ and $D^{\overline{F}}_z$ be the twisted Dirac operator along the fibre. Denote by ${\rm ch}(\overline{F})$ the Chern character of $\overline{F}$,  we can get a twisted
index
bundle and the twisted family index theorem:\\

\noindent{\bf Theorem 2.13}~(Atiyah-Singer [AS]) {\it It holds that the form}
$$(2i\pi)^{-\frac{n}{2}}\int_{M_z}\widehat{A}(R^{TZ}){\rm ch}(\overline{F})\eqno(2.39)$$
{\it is a representative of ${\rm {Ch}}([{\rm ker}D_{z,+}^{\overline{F}}]-[{\rm
ker}D_{z,-}^{\overline{F}}])$.}\\

\indent We shall prove Theorem 2.12 and Theorem 2.13
similarly as follows. For $t>0$, we set that $S_x$ is the spinor space on $T_x(M_z)$ and
$$\psi_t:~{\rm Hom}(S_x,S_x)\widehat{\otimes}\wedge(T_z^*B)\rightarrow {\rm Hom}(S_x,S_x)\widehat{\otimes}\wedge(T_z^*B)
~~hdy_\alpha\mapsto\frac{1}{\sqrt{t}}hdy_\alpha.\eqno(2.40)$$
For $\omega\in \wedge(T_z^*B)$ and $A\in {\rm Hom}(S_x,S_x)$, we define ${\rm Tr_s}(\omega A)=\omega {\rm Tr_s}(A).$
 By
Proposition 4.1 in [Z1], we get\\

\noindent{\bf Proposition 2.14} ([Bi2, Z1]) For $t > 0$, \\
$$\psi_t\int_{M_z}{\rm Tr_s}k_t^z(x,x)dx$$
{\it is a representative of ${\rm {Ch}}([{\rm ker}D_{z,+}]-[{\rm
ker}D_{z,-}])$.}\\

By Proposition 2.14, in order to prove Theorem 2.12, we only need to prove the following
local family index formula:\\

\noindent{\bf Theorem 2.15}([AS]) {\it We have}
$${\rm lim}_{t\rightarrow 0}\psi_t\int_{M_z}{\rm Tr}_sk_t^z(x,x)dx=(2i\pi)^{-\frac{n}{2}}\int_{M_z}\widehat{A}(R^{TZ}).\eqno(2.41)$$

In order to prove Theorem 2.15, let us recall the symbol map and the Getzler order, then we compute the symbol map of the heat kernel.
Let $n$ be even and ${\rm End}(S(TM_z))$ is a bundle of algebras
over $M_z$ isomorphic to the Clifford bundle Cl$(M_z)$, whose fiber
${\rm Cl}_x(M_z)$ at $x\in M_z$ is the complex algebra generated by
1 and elements of $T_x^*M_z$ with relations
$$\xi\cdot\eta+\eta\cdot\xi=-2\left<\xi,\eta\right>, ~~\xi,\eta\in
T_x^*M_z.\eqno(2.42)$$
 Recall that the quantization map $c:\wedge
 T^*_{{\mathbb{C}}}(M_z)\rightarrow {\rm Cl}(M)$ and the symbol map
 $\sigma=c^{-1}$ satisfy
 $$\sigma(c(\xi)c(\eta))=\xi\wedge\eta-\left<\xi,\eta\right>.\eqno(2.43)$$
So, for $\xi$ and $\eta$ in $\wedge T^*_{{\mathbb{C}}}(M_z)$ we have
$$\sigma(c(\xi^{(i)})c(\eta^{(j)}))=\xi^{(i)}\wedge \eta^{(j)}~~{\rm
mod}~~\wedge^{i+j-2} T^*_{{\mathbb{C}}}(M_z).\eqno(2.44)$$ where
$\xi^{(l)}$ denotes the component in $\wedge^{l}
T^*_{{\mathbb{C}}}(M_z)$ of $\xi\in \wedge T^*_{{\mathbb{C}}}(M_z).$
Since we compute the local index in a fixed fibre $M_z$, so we use an observation due to Getzler as in the case of a single manifold. If $e_1,\cdots,e_n$ is an orthonormal frame of $T_xM_z$,
then
$${\rm Tr}_s[c(e^{i_1})\cdots c(e^{i_k})]=
\left\{\begin{array}{lcr}
0 &~~{\rm if }~k\neq
n\\
(-2i)^{\frac{n}{2}}~&~{\rm if} ~~k=n.
\end{array}\right.
 \eqno(2.45)$$ We shall prove
(2.41) at a fixed point $x_0\in M_z$. Using normal coordinates
centered at $x_0$ in $M_z$ and paralleling $\partial_i$ at $x_0$
along geodesics through $x_0$, we get the orthonormal frame
$e_1,\cdots,e_n$. By
$$k_t^z(0,0)=K_Q^z(0,0,t)+O(t^{\infty})~~~~~~~~{\rm
as}~~t\rightarrow 0^+,$$ and by (2.45), we get
$${\rm lim}_{t\rightarrow 0}\psi_t\int_{M_z}{\rm Tr}_sk_t^z(x,x)dx=(-2i)^{\frac{n}{2}}{\rm lim}_{t\rightarrow 0}\int_{M_z}
\sigma[\psi_tK_Q(0,0,t)].\eqno(2.46)$$
 We define the Getzler order in [Po1] as follows:
 $${\rm deg}\partial_j=\frac{1}{2}{\rm deg}\partial_t={\rm deg}c(dx_j)={\rm deg}c(dy_j)=-{\rm deg}x^j=1.\eqno(2.47)$$
 Let $Q\in \Psi_V^*({\mathbb{R}}^n\times {\mathbb{R}}, S(TM_z)\otimes \wedge ^*T^*_zB)$ have symbol $$q(x,\xi,\tau)\sim
 \sum_{k\leq m'}\sum_{l=1}^{2^{\overline{q}}}q_{k,l}(x,\xi,\tau)\omega^{[l]},\eqno(2.48)$$
 where $\omega^{[l]}\in \wedge ^lT^*_zB$ and
$q_{k,l}(x,\xi,\tau)$ is an order $k$ symbol. Then taking components
in each subspace $\wedge^jT^*M_z$
and using Taylor expansions at $x = 0$ give formal expansions
$$\sigma[q(x,\xi,\tau)]\sim\sum_{j,k,l}\sigma[q_{k,l}(x,\xi,\tau)]^{(j)}\omega^{[l]}\sim\sum_{j,k,l,\alpha}\frac{x^\alpha}{\alpha!}
\sigma[\partial_x^\alpha
q_{k,l}(0,\xi,\tau)]^{(j)}\omega^{[l]}.\eqno(2.49)$$ The symbol
$\frac{x^\alpha}{\alpha!} \sigma[\partial_x^\alpha
q_{k,l}(0,\xi,\tau)]^{(j)}\omega^{[l]}$ is Getzler homogeneous
of order $k+j+l-|\alpha|$. So we can expand $\sigma[q(x,\xi,\tau)]$ as
$$\sigma[q(x,\xi,\tau)]\sim \sum_{j\geq 0}q_{(m-j)}(x,\xi,\tau),~~~~~~~~~q_{(m)}\neq 0, \eqno(2.50)$$
where $q_{(m-j)}$ is a Getzler homogeneous symbol of degree $m-j$.\\

\noindent {\bf Definition 2.16} ([Po1]) The integer $m$ is called the Getzler order of $Q$. The symbol $q_{(m)}$ is the principle Getzler
homogeneous symbol of $Q$. The operator $Q_{(m)}=q_{(m)}(x,D_x,D_t)$ is called the model operator of $Q$.\\

\noindent Denote by $O(t^{\frac{k}{2}})$ a Laurant expansion about $t^{\frac{1}{2}}$ whose lowest degree about $t$ is $\frac{k}{2}.$\\

\noindent {\bf Lemma 2.17} {\it Let $Q\in \Psi_V^*({\mathbb{R}}^n\times {\mathbb{R}}, S(T(M_z))\otimes \wedge ^*T^*_zB)$ have Getzler order $m$
and model operator $Q_{(m)}$. Then as $t\rightarrow 0^+$ we have:}
$$ \sigma[\psi_tK_Q(0,0,t)]^{(j)}=\omega^{{\rm odd}}O(t^{\frac{j-n-m-2}{2}})+O(t^{\frac{j-n-m-1}{2}}),~~{ if~} m-j~~{~is~ odd};$$
$$\sigma[\psi_tK_Q(0,0,t)]^{(j)}=t^{\frac{j-n-m-2}{2}}K_{Q_{(m)}}(0,0,1)^{(j)}+\omega^{{\rm odd}}_1O(t^{\frac{j-n-m-1}{2}})+O(t^{\frac{j-n-m}{2}}),$$
${ if~} m-j~~{~is~ even},$ {\it where $\sigma[K_Q(0,0,t)]^{(j)}$ denotes the degree $j$ form component in $M_z$ and $\omega^{{\rm odd}},$ $\omega^{{\rm odd}}_1$ are in $\wedge^{{\rm odd}}(T^*B)\otimes \wedge(T^*(M_z))$.
In particular $m=-2$ and $j=n$ is even, we get}
$$\sigma[\psi_tK_Q(0,0,t)]^{(n)}=K_{Q_{(-2)}}(0,0,1)^{(n)}+O(t).\eqno(2.51)$$\\
{\bf Proof.} By (1.7) in [Po1], we get
$$K_Q(0,0,t)\sim\sum_{m_0-j_0~{\rm even}}\sum_{l}t^{\frac{j_0-n-m_0-2}{2}}\check{q}_{m_0-j_0,l}(0,0,1)\omega^{[l]},\eqno(2.52)$$
where $m_0$ is the operator order of $Q$. And then
$$\sigma[\psi_tK_Q(0,0,t)]^{(j)}\sim\sum_{m_0-j_0~{\rm even}}\sum_{l}t^{\frac{j_0-n-m_0-l-2}{2}}\sigma[\check{q}_{m_0-j_0,l}(0,0,1)]^{(j)}\omega^{[l]}.\eqno(2.53)$$
Let $L=m_0-j_0+j+l$. By $Q$ having the Getzler order $m$, then
$L\leq m$. Then
$$\sigma[\psi_tK_Q(0,0,t)]^{(j)}\sim\sum_{m_0-j_0~{\rm even}}\sum_{l}\sum_{L\leq m}t^{\frac{j-n-L-2}{2}}\sigma[\check{q}_{m_0-j_0,l}(0,0,1)]^{(j)}\omega^{[l]},
\eqno(2.54)$$ We note that the leading term degree is $L=m$ and
$m_0-j_0+l=m-j$. When $m-j$ is odd, since $m_0-j_0$ is even, then
$l$ is odd, then we get
$$\sigma[\psi_tK_Q(0,0,t)]^{(j)}=\omega^{{\rm odd}}O(t^{\frac{j-n-m-2}{2}})+O(t^{\frac{j-n-m-1}{2}}).\eqno(2.55)$$
When $L=m$ and $m-j$ are even, $l$ is even. In this case, the leading coefficient is
$$\sigma[\breve{q}_{(m)}(0,0,1)]^{(j)}=\sum_{l}\sigma[\breve{q}_{m-j-l,l}(0,0,1)]^{(j)}\omega^{[l]}=K_{Q_{(m)}}(0,0,1)^{(j)}.\eqno(2.56)$$
In the next term, $L=m-1$ and $m-j$ are even and $m_0-j_0+l+j=m-1$, then $l$ is odd, so $m-j-l$ is odd and $\breve{q}_{m-j-l,l}(0,0,1)=0.$
Then the next term is $O(t^{\frac{j-n-m}{2}})$. ~~~~~~~~~~~~~~~~~~$\Box$

In the sequel we say that a symbol or a $\Psi DO$ is $O_G(m)$ if it has Getzler order $\leq m.$ Similar to Lemma 4 in [Po1], we have:\\

\noindent {\bf Lemma 2.18}~ {\it For $j = 1, 2$, let $Q_j\in
\Psi_V^*({\mathbb{R}}^n\times {\mathbb{R}}, S(T(M_z))\otimes \wedge
^*T^*_zB)$ have Getzler order $m_j$ and model operator $Q_{(m_j)}$
and assume that either $Q_1$ or $Q_2$ is properly supported. Then we have:}
$$ Q_1Q_2 = c[Q_{(m_1)}Q_{(m_2)}] + O_G(m_1 + m_2-1).\eqno(2.57)$$\\

\noindent{\bf Proof.} Let $Q_j=\omega_j\widehat{Q}_j$ for $j=1,2$ where $\omega_j\in \wedge
^{l_j}(T^*_zB)$ and $\widehat{Q}_j\in \Psi_V^*({\mathbb{R}}^n\times {\mathbb{R}}, S(T(M_z))$ has the Getzler order $m_j-l_j$. Then
$$ Q_1Q_2 =\omega_1\wedge\omega_2\widehat{Q}_1\widehat{Q}_2.\eqno(2.58)$$
By Lemma 4 in [Po1], we have
$$\widehat{Q}_1\widehat{Q}_2=c[\widehat{Q}_{(m_1-l_1)}\widehat{Q}_{(m_2-l_2)}]+O_G(m_1+m_2-l_1-l_2-1).\eqno(2.59)$$
By (2.58) and (2.59), we get (2.57).~~~~~~~~~~~$\Box$\\

\indent In the following, we compute the model operator of $F$ in (2.5). Let $x_0,x\in M_z$ and $\tau^E(x_0,x)$ be the parallel transport map
in the bundle $\pi^*\wedge^*(T_z^*B)\otimes S(TM_z)$ along the geodesic from $x$ to $x_0$, defined with respect to the Clifford connection
$\nabla^{E,0}$. Using this map, we can trivialize the bundle $\pi^*\wedge^*(T_z^*B)\otimes S(TM_z)$. We note that $e_1,\cdots, e_n$ is
the parallel transport frame with respect to $\nabla^{TZ}$ and so $m_0(e_i)$ is not a constant matrix under above trivialization. But we have\\

\noindent {\bf Lemma 2.19} ([BGV Lemma 10.25]) {\it If $c^i$ is the Clifford action of the cotangent vector $dx_i$ acting on $E$ at $x_0$ and
$\varepsilon^\alpha$ is multiplication by $f^\alpha$ in the exterior algebra $\wedge^*(T_z^*B)$, we have}
$$m_0(e^i)=c^i+\sum u_\alpha^i\varepsilon^\alpha;~~~~~m_0(f^\alpha)=\varepsilon^\alpha,\eqno(2.60)$$
{\it where $u_\alpha^i$ are smooth functions on $U$ satisfying $u_\alpha^i(x)=O(|x|).$}\\

\noindent {\bf Lemma 2.20} ([BGV Lemma 10.26]) {\it In the trivialization of $E$ over $U$ induced by the parallel  transport map $\tau^E(x_0,x)$,
the connection $\nabla^{E,0}$ equals $d+\Theta$, where}
$$\Theta(\partial_i)=-\frac{1}{4}\sum_{j,a<b}\left<R(\partial_i,\partial_j)e_a,e_b \right> m^am^bx^j+\sum_{a<b}f_{iab}(x)m^am^b+g_i(x);\eqno(2.61)$$
{\it here $m^a$ represents $c^i$ or $\varepsilon^\alpha$ and $f_{iab}=O(|x|^2)$ and $g_i(x)=O(|x|)$.}\\

By (2.5) and Lemma 2.19 and Lemma 2.20, we get \\

\noindent {\bf Proposition 2.21} {\it In the trivialization of $E$ over $U$ induced by the parallel transport map $\tau^E(x_0,x)$ and the
normal coordinates, the model operator of $F$ is}
$$F_{(2)}=-\sum_{i=1}^n(\partial_i-\frac{1}{4}\sum_{j=1}^na_{ij}x_j)^2,~~~~a_{ij}=\left<R^{TZ}\partial_i,\partial_j\right>.\eqno(2.62)$$\\

We note that ${\rm lim}_{t\rightarrow 0}\psi_t\int_{M_z}{\rm Tr}_sk_t^z(x,x)dx$ does not depend on the trivialization. Using the same discussions with Lemma 5 and Lemma 6 in [Po1], we get\\

\noindent {\bf Lemma 2.22} {\it Let $Q$ be a Volterra paramatrix for $F+\partial_t$. Then $Q$ has the model operator $(F_{(2)}+\partial_t)^{-1}$ and}
$$K_{(F_{(2)}+\partial_t)^{-1}}(0,0,1)=(4\pi)^{-\frac{n}{2}}\widehat{A}(R^{TZ}).\eqno(2.63)$$\\

By (2.46), (2.51) and Lemma 2.22, we proved Theorem 2.15.\\

\section { The local equivariant family index theorem}

\quad In [LM], Liu and Ma gave a proof of the local equivariant family index theorem by the finite propagation speed method in [B3]. In [PW],
 a new proof of the local equivariant index theorem was given by the Volterra calculus. In this section, we shall give a new proof of the local equivariant family index theorem by the Volterra calculus.\\
 \indent Let us give the fundamental setup. We assume that $M$ and $TZ$ are oriented and $G$ is a compact Lie group which is a fiberwise isometry on $M$ and preserves the orientation of $TZ$. Then $G$ acts as identity on $B$. We also assume that the action of G lifts to $S(TZ)$
and that the $G$-action commutes with $D$. By Proposition 1.1 in
[LM], we know that ${\rm Ind} D^G\in K_G(B)$.  Now let us calculate
the equivariant Chern character $\rm {ch}_\phi({\rm Ind}(D^G))$ in
terms of the fixed point data of $\phi\in G$. Set $ M^\phi =\{x\in
M, \phi x = x\}$. Then $\pi: M^\phi \rightarrow B$ is a fibration
with compact fibre $M^\phi_z$. By [BGV, Proposition 6.14], $TZ^\phi$
is naturally oriented in $M^\phi$. Let $N$ denote the normal bundle
of $M^\phi$, then $N =TZ/TZ^\phi$. We fixed a fibre $M_z$ and denote
by $\widetilde{\phi}$ the lift of $\phi$ which maps $S(T(M_z)_x)$ to
$S(T(M_z)_{\phi x})$.
   We denote by $M^\phi_z$ the fixed-point set of $\phi$, and for $a = 0,\cdots ,n,$ we let
   $M^\phi_z=\bigcup _{0\leq a\leq n} M_{z,a}^\phi$, where  $M_{z,a}^\phi$ is an $a$-dimensional submanifold. Given a fixed-point $x_0$ in a component
   $M_{z,a}^\phi$, consider some local coordinates $x = (x^1,\cdots , x^a)$ around
$x_0.$ Setting $b = n-a,$ we may further assume that over the range
of the domain of the local coordinates there is an orthonormal frame
$e_1(x),\cdots , e_b(x)$ of $N^\phi_z$. This defines fiber
coordinates $v = (v_1, \cdots , v_b).$ Composing with the map
$(x,v)\in N^\phi_z(\varepsilon_0)\rightarrow {\rm exp}_x(v)$ we then
get local coordinates $x^1,\cdots,x^a,v^1,\cdots,v^b$ for $M_z$ near
the fixed point $x_0$. We shall refer to this type of coordinates as
{\it tubular coordinates.} Then $N^\phi_z(\varepsilon_0)$ is homeomorphic with a tubular neighborhood of $M^\phi_z$. We have \\

\noindent{\bf Theorem 3.1} ([LM]) {\it For any $t > 0$, the form
${\rm Str}[\widetilde{\phi}\psi_t{\rm exp}(-tF)]$ is closed and
its de-Rham cohomology
class in $B$ is independent of $t$ and represents ${\rm ch}_\phi({\rm Ind}(D^G))$ in the de-Rham cohomology of $B$.}\\

\indent We shall use Theorem 3.1 to find a local index formula for ${\rm ch}_\phi({\rm Ind}(D^G))$ by estimating ${\rm Str}[\widetilde{\phi}\psi_t{\rm exp}(-tF)]$. By the Mckean-Singer formula,
we have
$${\rm Str}[\widetilde{\phi}\psi_t{\rm exp}(-tF)]=\psi_t\int_{M_z}{\rm Str}[ \widetilde{\phi}k_t^z(x,\phi(x))]dx
=\psi_t\int_{M_z}{\rm Str}[
\widetilde{\phi}K_{(F+\partial_t)^{-1}}^z(x,\phi(x),t)]dx.
\eqno(3.1)$$
 Let $Q=(F+\partial_t)^{-2}$. For $x\in M_z^\phi$ and
$t>0$ set
$$I_Q(x,t):=\widetilde{\phi}(x)^{-1}\int_{N_x^\phi(\varepsilon)}\widetilde{\phi}({\rm
exp}_xv)K_Q({\rm exp}_xv,{\rm exp}_x(\phi'(x)v),t)dv.\eqno(3.2)$$
Here we use the trivialization of $S(TM_z)$ about the tubular
coordinates. Using the tubular coordinates, then
$$I_Q(x,t)=\int_{|v|<\varepsilon}\widetilde{\phi}(x,0)^{-1}\widetilde{\phi}(x,v)K_Q(x,v;x,\phi'(x)v;t)dv.\eqno(3.3)$$
Let
$$q^E_{m-j}(x,v;\xi,\nu;\tau):=\widetilde{\phi}(x,0)^{-1}\widetilde{\phi}(x,v)q_{m-j}(x,v;\xi,\nu;\tau).\eqno(3.4)$$
Using the same proof as Lemma 2.9, we obtain the following proposition by Lemma 8.2 and Lemma 8.3 in [PW]\\

\noindent {\bf Proposition 3.2} {\it Let $Q\in
\Psi_V^m(M_z\times {\mathbb{R}},E),~m\in {\mathbb{Z}}.$ Uniformly on
each component $M_{z,a}^\phi$
$$I_Q(x,t)\sim \sum_{j \geq
0}t^{-(\frac{a}{2}+[\frac{m}{2}]+1)}I_Q^j(x) ~~~~~~{\rm
as}~~t\rightarrow 0^+,\eqno(3.5)$$ where $I_Q^j(x)$ is defined in tubular coordinates by}
$$I_Q^{(j)}(x):=\sum_{|\alpha|\leq
m-[\frac{m}{2}]+2j}\int\frac{v^\alpha}{\alpha!}\left(\partial_v^\alpha
q^E_{2[\frac{m}{2}]-2j+|\alpha|}\right)^\vee(x,0;0,(1-\phi'(x))v;1)dv.\eqno(3.6)$$\\

\indent By Proposition 8.7 in [PW] and Proposition 3.2 and the definition of $\psi_t$, we have\\

\noindent {\bf Proposition 3.3} {\it  Let $P:~C^{\infty}(M_z,
\pi^*\wedge T^*_zB\otimes S(T(M_z)))\rightarrow C^{\infty}(M_z,
\pi^*\wedge T^*_zB\otimes S(T(M_z)))$ be a differential operator of
order
$m$.}\\
{\rm (1)}{\it  Uniformly on each component $M_{z,a}^\phi$,}
$$I_{P(F+\partial_t)^{-1}}(x,t)\sim\sum_{j\geq 0}
t^{-(\frac{n}{2}+[\frac{m}{2}])+j}I^{(j)}_{P(F+\partial_t)^{-1}}(x)~~~~~{\rm
as}~~t\rightarrow 0^+,\eqno(3.7)$$
 {\rm (2)}{\it  As $t\rightarrow 0^+$, we have}
 $$\psi_t{\rm Str}[\widetilde{\phi}Pe^{-tF}]\sim\sum_{0\leq a\leq
 n}\sum_{j \geq
 0}\sum_{l=1}^{2^{\overline{q}}}t^{-(\frac{a}{2}+[\frac{m}{2}]+\frac{l}{2})+j}\int_{M_{z,a}^\phi}
 {\rm
Str}[\widetilde{\phi}I^{(j),l}_{P(F+\partial_t)^{-1}}(x)]dx,\eqno(3.8)$$
{\it where $I^{(j),l}_{P(F+\partial_t)^{-1}}(x)$ denotes the
$l$-degree component of $I^{(j)}_{P(F+\partial_t)^{-1}}(x)$ in
$\wedge
T^*_zB.$}\\

Let $$\widehat{A}(R^{TZ^\phi})={\rm
det}^{\frac{1}{2}}\left(\frac{R^{TZ^\phi}/2}{{\rm
sinh}(R^{TZ^\phi}/2)}\right);~~\nu_\phi(R^{N^\phi}):={\rm
det}^{-\frac{1}{2}}(1-\phi^Ne^{-R^{N^\phi}}).\eqno(3.9)$$
For a top degree form $\omega \in C^{\infty}(M_{z,a},\wedge T ^*_z(B) \otimes \wedge ^a T^*M_{z,a})$ in $M_{z,a}$ with coefficients in $\wedge T^*_z(B)$,
we denote by $|\omega|^{(a)}$ the Berezin integral which in $\wedge T ^*_z(B)$, i.e. its inner-product with the volume form of $M_{z,a}.$\\

\noindent {\bf Theorem 3.4} (Local equivariant family index theorem)
{\it Let $x_0\in M^\phi$, then}\\
$${\rm lim}_{t\rightarrow 0}{\rm
Str}\left[\widetilde{\phi}(x_0)\psi_tI_{(F+\partial_t)^{-1}}(x_0,t)\right]=(-i)^{\frac{n}{2}}(2\pi)^{-\frac{a}{2}}\left|\widehat{A}(R^{TZ^\phi})
\nu_\phi(R^{N^\phi})\right|^{(a)}.\eqno(3.10)$$ \\

\noindent {\bf Remark.} It is easy to generalize Theorem 3.4 to the twisted $\phi$-complex vector bundle case.\\

\indent By (3.1) and Lemma 8.1 in [PW], as $t\rightarrow 0^+$
$$\psi_t\int_{M_z}{\rm Str}[
\widetilde{\phi}K_{(F+\partial_t)^{-1}}(x,\phi(x),t)]dx
=\psi_t\int_{M_z}{\rm Str}[
\widetilde{\phi}I_{Q}(x,\phi(x),t)]dx,\eqno(3.11)$$
 By Theorem 3.1, Theorem 3.4 and (3.1), (3.11), we get a representative of ${\rm ch}_\phi({\rm Ind}D^G)$.\\
\indent In order to prove Theorem 3.4, we shall compute $\widetilde{\phi}$ in tubular coordinates.
 Let $e_1, \dots , e_n$ be an oriented orthonormal basis of
$T_{x_0}M_z$ such that $e_1,\cdots , e_a$ span $T_{x_0}M_z^\phi$ and
$e_{a+1},\cdots , e_n$ span $N_{ x_0}^\phi$ . This provides us with
normal coordinates $(x_1, \cdots , x_n)\rightarrow {\rm
exp}_{x_0}(x^1e_1+\cdots+x^ne_n).$ Moreover using parallel
translation enables us to construct a synchronous local oriented
tangent frame $e_1(x), . . . , e_n(x)$ such that $e_1(x),\cdots ,
e_a(x)$ form an oriented frame of $TM_{z,a}^\phi$ and $e_{a+1}(x),
\cdots , e_n(x)$ form an (oriented) frame $N^\phi$ (when both frames
are restricted to $M_z^\phi).$ This gives rise to trivializations of
the tangent and spinor bundles. Write
$$\phi'(0)=\left(\begin{array}{lcr}
  1  & 0  \\
   0  &  \phi^N
\end{array}\right).
$$
Let $\wedge(n)=\wedge_{\mathbb{C}}^*{\mathbb{R}}^n$ be the
complexified exterior algebra of ${\mathbb{R}}^n$.We shall use the
following gradings on $\wedge(n),$
$$\wedge(n)=\bigoplus_{1\leq j\leq n}\wedge^j(n)=\bigoplus_{\begin{array}{lcr}
  1\leq k\leq a \\
  1\leq \overline{l} \leq b
\end{array}}\wedge^{k,\overline{l}}(n),$$
 where $\wedge^j(n)$ is the space of forms of degree $j$ and
 $\wedge^{k,\overline{l}}(n)$
is the space of forms $dx^{i_1}\wedge\cdots\wedge
dx^{i_{k+\overline{l}}}$ with $1\leq i_1<\cdots <i_k\leq a$ and $a +
1\leq i_{k+1} < \cdots < i_{k+\overline{l}}\leq n.$ Given a form
$\omega\in\wedge (n)$ we shall denote by $\omega^{(j)} $ (resp.,
$\omega^{(k,\overline{l})}$) its component in $\wedge^j(n) $
(resp.,$\wedge^{k,\overline{l}}(n)$ ). We denote by
$|\omega|^{(a,0)}$ the Berezin integral $|\omega^{(*,0)}|^{(a,0)}$
of its component $\omega^{(*,0)}$ in $\wedge^{(*,0)}(n).$ Then we have\\

\noindent{\bf Lemma 3.5 }([PW,Lemma 9.5]) {\it Let $S_n$ be the spinor space associated to ${\mathbb{R}}^n$ and $A\in {\rm End}(S_n)$, then}
$${\rm
Str}[\widetilde{\phi}A]=(-2i)^{\frac{n}{2}}2^{-\frac{b}{2}}{\rm
det}^{\frac{1}{2}}(1-\phi^N)|\sigma(A)|^{(a,0)}+(-2i)^{\frac{n}{2}}\sum_{0\leq
b'<b}
|\sigma(\widetilde{\phi})^{(0,b')}\sigma(A)^{(a,b-b')}|^{(n)}.\eqno(3.12)$$\\

Similar to Proposition 2.21, we get the same expression of the model
operator of $F$ in the trivialization of $E$ over $U$ induced by the
parallel transport map $\tau(x_0,x)$ about $\nabla^{TZ}$ and the
normal coordinate. We will compute the local index in this
trivialization.\\

\noindent{\bf Lemma 3.6}~ {\it Let the operator $Q\in \Psi_V^*({\mathbb{R}}^n\times
{\mathbb{R}}, S)\otimes \wedge(T_z^*B)$ have the Getzler order $m$
and model operator $Q_{(m)}$. Then as $t\rightarrow 0^+$}\\
\indent (1) $\sigma[\psi_tI_Q(0,t)]^{(j)}=\omega^{{\rm
odd}}O(t^{\frac{j-m-a-2}{2}})+O(t^{\frac{j-m-a-1}{2}})$
{\it if
$m-j$
is odd.}\\
\indent (2)
$\sigma[\psi_tI_Q(0,t)]^{(j)}=O(t^{\frac{j-m-a-2}{2}})I_{Q(m)}(0,1)^{(j)}+\omega_1^{{\rm
odd}}O(t^{\frac{j-m-a-1}{2}})+O(t^{\frac{j-m-a}{2}})$\\
\indent {\it ~~~~if
$m-j$
is even.}\\
{\it In particular, for $m=-2$ and $j=a$ we get}
$$\sigma[\psi_tI_Q(0,t)]^{(a,0)}=I_{Q(-2)}(0,1)^{(a,0)}+O(t^{\frac{1}{2}}).\eqno(3.13)$$

\noindent{\bf Proof.} By the change of variable formula for symbol,
similar to (9.23) in [PW], we have
$$\sigma[\psi_tI_Q(0,t)]^{(j)}\sim\sum_{\begin{array}{lcr}
 |\alpha|+|\beta|-|\gamma|-\overline{l} ~{\rm even}\\
  \overline{1}\leq m',~ 2|\gamma|<|\beta|
\end{array}}\sum_{l=1}^{2^{\overline{q}}}t^{\frac{|\alpha|+|\beta|-|\gamma|-l-\overline{l}-(a+2)}{2}}I_{\overline{l}\alpha\beta\gamma,l}^{(j)},\eqno(3.14)$$
where $m'$ is the operator order of $Q$ and
$$I_{\overline{l}\alpha\beta\gamma,l}^{(j)}=\int_{{\mathbb{R}}^b}a_{\beta\gamma,\overline{l}}(0,v)\frac{v^{\alpha}}{\alpha!}\left(
\partial_v^\alpha\sigma[\xi^\gamma D_\xi^\beta
q_{\overline{l},l}]^{(j)}\right)^\vee(0,0;0,(1-\phi^N(0))v;1)dv\omega^{[l]},\eqno(3.15)$$
where $a_{\beta\gamma,\overline{l}}(x',v)$ are smooth functions
such that $a_{\beta\gamma,\overline{l}}(x)=1$ when $\beta=\gamma=0.$
Since $Q$ has the Getzler order $m$, then all the coefficients
$I_{\overline{l}\alpha\beta\gamma,l}^{(j)}$ with
$\overline{l}+j+l-|\alpha|>m$ are zero. So, if
$\overline{l}+j+l-|\alpha|\leq m$ and $2|\gamma|\leq |\beta|,$ then
$t^{\frac{|\alpha|+|\beta|-|\gamma|-l-\overline{l}-(a+2)}{2}}$ is
$O(t^{\frac{j-m-a-2}{2}})$ and even is $o(t^{\frac{j-m-a-2}{2}})$ if
we have $\overline{l}+j+l-|\alpha|<m$ or $(\beta,\gamma)\neq
(0,0).$\\
\indent Observe that terms in asymptotic (3.14) containing integer
powers (resp. half integer powers) of $t$ when $l$ is even (resp.
odd) since $a$ is even. When $m-j$ is odd, the leading term is
$O(t^{\frac{j-m-a-2}{2}})$ which is a half integer, so its
coefficient is in $\Omega^{\rm odd}(B)$ and (1) is verified. When $m-j$
is even, similar discussions show that
$$\sigma[\psi_tI_Q(0,t)]^{(j)}=\omega^{{\rm
even}}O(t^{\frac{j-m-a-2}{2}})+\omega^{{\rm
odd}}O(t^{\frac{j-m-a-1}{2}})+O(t^{\frac{j-m-a}{2}}).\eqno(3.16)$$
Then similar to the discussions in [PW], we prove (2). $\Box$\\

By (3.12) and Lemma 3.6 (2) and (3.13) and Lemma 9.13 in [PW], we
get\\ Theorem 3.4.

\section{ The equivariant JLO character for a family of the Dirac
operators}
  \quad In [CH], Chern and Hu computed the equivariant JLO characters for
  invariant Dirac operators. In [Az2], Azmi computed the JLO characters for
 a family of Dirac operators. In [Az3], Azmi constructed an equivariant bivariant cyclic
theory, as a combination of equivariant cyclic and noncommutative de
Rham theories for unital $G$-Banach algebras, where $G$ is a compact
Lie group. By incorporating the JLO formula and the superconnection
formalism of Quillen, an equivariant bivariant JLO character of
Kasparov's G-bimodule is defined, with values in the bivariant
cyclic theory. In this section, following the ideas in [CH], [Zh2]
and [Fe], we compute the equivariant JLO characters for a family of
the
Dirac operators.\\
  \indent Let $C^1(M)$ be the Banach algebra which is the completion of $C^{\infty}(M)$ with respect
to the norm $|f| := ||f||+||[D,f]||,$ for $f\in C^{\infty}(M)$. The
commutator $[D, f]$ extends to a bounded operator on
${\mathcal{H}}=L^2(M,S(TZ))$. The algebra $C^1(M)$ acts on
$L({\mathcal{H}})$ \\(bounded operators on ${\mathcal{H}}$) by
multiplication.  Denote by $\Delta$ the
projective tensor product of the Banach algebras $C^1(M)$ and
$C^\infty(B)$, i.e $\Delta= C^\infty(B) \hat{\otimes}C^1(M),$ with
the
projective tensor product norm.\\
\indent Let ${\mathcal{M}}=\wedge B\otimes {\mathcal{H}}$  be a
$\Delta-C^{\infty}(B)$ bimodule, where $\Delta$ acts on the left of
${\mathcal{M}}$ by letting $C^\infty(B)$ act on $\wedge(B)$ by
multiplication by left and $C^1(M)$ acts on ${\mathcal{H}}$ while
$C^{\infty}(B)$ acts on ${\mathcal{M}}$ by multiplication from right.
There is a continuous $C^\infty(B)$-valued inner product on
${\mathcal{M}}.$ There is an obvious continuous action of $\phi\in
G$ on ${\mathcal{M}}$, by letting $\phi$ act on ${\mathcal{H}}$ via
$\widetilde{\phi}$ and on $\wedge B$ via the identity map. Let
${\mathcal{B}}_t=\sqrt{t}\psi_t({\mathcal{B}})$ be the rescaled
Bismut' superconnection with the rescaled curvature
$F_t={\mathcal{B}}_t^2=t\psi_t(F).$
If $h_i\otimes f_i\in \Delta,~ 0\leq i\leq 2k$ are operators on ${\mathcal{M}}$, we define the equivariant bivariant JLO character by:\\
$${\rm Ch}_{2k}({\mathcal{B}}_t)(\phi)(h_0\otimes f_0,\cdots,h_{2k}\otimes f_{2k})
=t^k\int_{\triangle_{2k}}{\rm
Str}\left[\psi_t\widetilde{\phi}h_0\otimes f_0e^{-ts_1F}\right.$$
$$
\left.\cdot [{\mathcal{B}}, h_1\otimes f_1] e^{-t(s_2-s_1)F}\cdots
[{\mathcal{B}}, h_{2k}\otimes f_{2k}]
e^{-t(1-s_{2k})F}\right]ds,\eqno(4.1)$$ where
$\triangle_{2k}=\{(s_1,\cdots,s_{2k})|~0\leq s_1\leq\cdots\leq
s_{2k}\leq 1\}$ is the simplex in ${\mathbb{ R}}^{2k}$. We will
compute $${\rm lim}_{t\rightarrow 0}{\rm
Ch}_{2k}({\mathcal{B}}_t)(\phi).$$ \indent In the following, we give
some estimates about ${\rm Ch}_{2k}({\mathcal{B}}_t)(\phi).$ Since
$B$ is compact, we can fix a fibre $M_z$ and estimate ${\rm
Ch}_{2k}({\mathcal{B}}_t)(\phi).$ Let $H$ be a Hilbert space. For
$q\geq0$, denote by $||.||_q$ the Schatten $p$-norm on the Schatten
ideal
$L^p$. Denote by $L(H)$ the Banach algebra of bounded operators on $H$.\\

\noindent {\bf Lemma 4.1}~([Si]){\it~~(i)~~${\rm Tr}(AB)={\rm
Tr}(BA)$, for $A,~B\in L(H)$ and $AB, ~BA\in
L^1$.\\
~~(ii)~~For $A\in L^1,$ we have $|{\rm Tr}(A)|\leq ||A||_1$,
$||A||\leq ||A||_1$.\\
~~(iii)~~For $A\in L^q$ and $B\in L(H)$, we have: $||AB||_q\leq
||B||||A||_q$, $||BA||_q\leq ||B||||A||_q$.\\
 ~~(iv)~(H\"{o}lder
Inequality)~~If $\frac{1}{r}=\frac{1}{p}+\frac{1}{q},~p,q,r>0,~A\in
L^p,~B\in
L^q,$ then $AB\in L^r$ and $||AB||_r\leq ||A||_p||B||_q$.}\\

\noindent{\bf Lemma 4.2} {\it For any $u>0,~t>0$ and $t$ is small
and any order $l$ fibrewise
differential operator $\overline{B}$ with form coefficients, we have:}\\
$$||e^{-utF}\overline{B}||_{u^{-1}}\leq C_lu^{-\frac{l}{2}}t^{-\frac{l}{2}}({\rm
tr}[e^{-\frac{tD^2}{2}}])^u.\eqno(4.2)$$\\

\noindent {\bf Proof.} By (2.7), we have
$$||e^{-utF}\overline{B}||_{u^{-1}}=||
\sum_{m\geq
0}(-ut)^m\int_{\triangle_m}e^{-v_0utD^2}F_{[+]}e^{-v_1utD^2}$$
$$\cdot F_{[+]}\cdots e^{-v_{m-1}utD^2}F_{[+]}e^{-v_mutD^2}\overline{B}dv||_{u^{-1}},\eqno(4.3)$$
We estimate the term of $m=2$ in the right hand of (4.3), other
terms are similar. We split $\triangle(2)=J_0\cup J_1 \cup J_2$
where $J_i=\{(v_0,v_1,v_2)\in \Delta(2)|v_i\geq \frac{1}{3}\}.$
\begin{eqnarray*}
&&(ut)^2||\int_{J_0}e^{-v_0utD^2}F_{[+]}e^{-v_1utD^2}
F_{[+]}e^{-v_2utD^2}\overline{B}dv||_{u^{-1}}\\
&&\leq (ut)^2\int_{J_0}||e^{-\frac{v_0ut}{2}D^2}||_{(uv_0)^{-1}}
||e^{-\frac{v_0ut}{2}D^2}(1+D^2)^{\frac{l+2}{2}}||||(1+D^2)^{-\frac{l+2}{2}}
F_{[+]}(1+D^2)^{\frac{l+1}{2}}||\\
&&\cdot||e^{-{v_1ut}D^2}||_{(uv_1)^{-1}}||(1+D^2)^{-\frac{l+1}{2}}
F_{[+]}(1+D^2)^{\frac{l}{2}}||||e^{-{v_2ut}D^2}||_{(uv_2)^{-1}}||1+D^2)^{-\frac{l}{2}}\overline{B}||dv\\
&&\leq (ut)^2\int_{J_0} \left({\rm
Tr}e^{-\frac{t}{2}D^2}\right)^{uv_0}\left({\rm
Tr}e^{-{t}D^2}\right)^{u(v_1+v_2)}(uv_0t)^{-\frac{l+2}{2}}\\
&&\cdot||(1+D^2)^{-\frac{l+2}{2}}
F_{[+]}(1+D^2)^{\frac{l+1}{2}}||||(1+D^2)^{-\frac{l+1}{2}}
F_{[+]}(1+D^2)^{\frac{l}{2}}||||1+D^2)^{-\frac{l}{2}}\overline{B}||dv\\
&&\leq
 C_2\left({\rm
Tr}e^{-\frac{t}{2}D^2}\right)^{u}(ut)^{-\frac{l}{2}+1}~~~~~~~~~~(4.4)
\end{eqnarray*}
where we use $F_{[+]}$ is a fibrewise first order differential
operator and the equality
$${\rm
sup}\{(1+x)^{\frac{l}{2}}e^{-\frac{utx}{2}}\}=(ut)^{-\frac{l}{2}}e^{-\frac{l-ut}{2}}.\eqno(4.5)$$
In the above estimate, we omit the norm of coefficient forms since
$B$ is compact. For $J_1$ and $J_2$ we have similar estimates. For
the general $m$, we get
$$||(-ut)^m\int_{\triangle_m}e^{-v_0utD^2}F_{[+]}e^{-v_1utD^2}F_{[+]}\cdots
e^{-v_{m-1}utD^2}$$
$$\cdot F_{[+]}e^{-v_mutD^2}\overline{B}dv||_{u^{-1}}\leq C_2\left({\rm
Tr}e^{-\frac{t}{2}D^2}\right)^{u}(ut)^{-\frac{l}{2}+\frac{m}{2}}.\eqno(4.6)$$
By (4.3) and (4.6), we get (4.2). $\Box$\\

 \noindent{\bf Lemma 4.3}{\it ~Let $\overline{B}_1,~\overline{B}_2$ be positive order
 $p,~q$ fibrewise pseudodifferential operators with form coefficients respectively, then for any
 $s,~t>0,~0\leq u\leq1$, we have the following estimate:}\\
$$||\overline{B}_1e^{-ustF}\overline{B}_2e^{-(1-u)stF}||_{s^{-1}}\leq C_{p,q}s^{-\frac{p+q}{2}}t^{-\frac{p+q}{2}}({\rm
tr}[e^{-\frac{tD^2}{4}}])^s.\eqno(4.7)$$\\

\noindent{\bf Proof.} Similar to the proof of Lemma 4.2, we have when $t$ is small and $u\geq \frac{1}{2}$,
$$||(1+D^2)^{-\frac{q}{2}}e^{-(u-\frac{1}{2})stF}(1+D^2)^{\frac{q}{2}}||\leq C_0;~~||e^{-(1-u)stF}||\leq C_1,\eqno(4.8)$$
where $C_0,C_1$ are constants which do not depend on $t,u.$ When $u\geq \frac{1}{2}$, by Lemma 4.2 and (4.8), we have:
\begin{eqnarray*}
&&||\overline{B}_1e^{-ustF}\overline{B}_2e^{-(1-u)stF}||_{s^{-1}}\\
&&\leq ||\overline{B}_1(1+D^2)^{-\frac{p}{2}}||||(1+D^2)^{\frac{p}{2}}e^{-\frac{1}{2}stF}(1+D^2)^{\frac{q}{2}}||_{s^{-1}}\\
&&||(1+D^2)^{-\frac{q}{2}}e^{-(u-\frac{1}{2})stF}(1+D^2)^{\frac{q}{2}}||||(1+D^2)^{-\frac{q}{2}}\overline{B}_2||||e^{-(1-u)stF}||\\
&&\leq C_{p,q}s^{-\frac{p+q}{2}}t^{-\frac{p+q}{2}}({\rm
tr}[e^{-\frac{tD^2}{4}}])^s,~~~~~~~~~~~~~~~~~~~~~~~~~~~~~~~~~~~~~~~~~~~~~~~~~~~~~~~~~~~~(4.9)
\end{eqnarray*}
where we use
$$|(1+D^2)^{\frac{p}{2}}e^{-\frac{1}{2}stF}(1+D^2)^{\frac{q}{2}}||_{s^{-1}}=||e^{-\frac{1}{2}stF}(1+D^2)^{\frac{p+q}{2}}||_{s^{-1}}.\eqno(4.10)$$
When $u\leq \frac{1}{2}$, by Lemma 4.2, we have:
\begin{eqnarray*}
&&||\overline{B}_1e^{-ustF}\overline{B}_2e^{-(1-u)stF}||_{s^{-1}}\\
&&\leq ||\overline{B}_1e^{-ustF}(1+D^2)^{-\frac{p}{2}}||||(1+D^2)^{\frac{p}{2}}\overline{B}_2e^{-\frac{1}{2}stF}||_{s^{-1}}
||e^{-(1-u-\frac{1}{2})stF}||\\
&&\leq C_{p,q}s^{-\frac{p+q}{2}}t^{-\frac{p+q}{2}}({\rm tr}[e^{-\frac{tD^2}{4}}])^s,~~~~~~~~~~~~~~~~~~~~~~~~~~~~~~~~~~~~~~~~~~~~~~~~~~~~~~~~~~~~(4.11)
\end{eqnarray*}
By (4.9) and (4.11), we have proved this lemma. $\Box$\\

\indent Let $\overline{B}$ be a fibrewise
operator with form coefficients and $l$ be a positive interger.
 Write
$$\overline{B}^{[l]}=[F,\overline{B}^{[l-1]}],~\overline{B}^{[0]}=\overline{B}.$$\\

\noindent{\bf Lemma 4.4}~{\it Let $\overline{B}$ a finite order fibrewise
differential
operator with form coefficients, then for any $s>0$, we have:}\\
$$e^{-sF}\overline{B}=\sum^{N-1}_{l=0}\frac{(-1)^l}{l!}s^l\overline{B}^{[l]}e^{-sF}+(-1)^Ns^N\overline{B}^{[N]}(s),\eqno(4.12)$$
\noindent {\it where $\overline{B}^{[N]}(s)$ is given by}\\
$$\overline{B}^{[N]}(s)=\int_{\triangle_N}e^{-u_1sF}\overline{B}^{[N]}e^{-(1-u_1)sF}du_1du_2\cdots du_N.\eqno(4.13)$$
\\

\noindent{\bf Proof.} Here we use Lemma 1.9 in [BeC],
$$[A, e^{-F}]=-\int_0^1e^{-sF}[A,F]e^{-(1-s)F}ds,\eqno(4.14)$$
then similar to Lemma 3 in [Fe], we prove this lemma. $\Box$\\

For $h_j\otimes f_j$, we have
$$[{\mathcal{B}},h_j\otimes f_j]=d_B(h_j\otimes f_j)\wedge+h\otimes [D,f].\eqno(4.15)$$
Write $T_j$ is $d_B(h_j\otimes f_j)\wedge$ or $h\otimes [D,f], 1\leq
j \leq 2k$ and $T_0=h_0\otimes f_0.$
By Lemma 4.4, we have:\\
$$T_0e^{-s_1tF}T_1e^{-(s_2-s_1)tF}T_2\cdots
e^{-(s_{2k}-s_{2k-1})tF}T_{2k}e^{-(1-s_{2k})tF}$$
$$=\sum^{N-1}_{\lambda _{1},\cdots ,\lambda _{2k}=0}\frac{(-1)^{\lambda _1+\cdots +\lambda _{2k}}{s_1}^{\lambda _{1}}\cdots
 s_{2k}^{\lambda
_{2k}}t^{\lambda _{1}+\cdots +\lambda _{2k}}} {\lambda
_{1}!\cdots\lambda _{2k}!} T_0[T_1]^{[\lambda _{1}]}\cdots
[T_{2k}]^{[\lambda _{2k}]}e^{-tF}$$
$$+\sum_{1\leq q\leq {2k}}\sum^{N-1}_{\lambda _{1},\cdots ,\lambda _{q-1}=0}
\frac{(-1)^{\lambda _1+\cdots +\lambda _{q-1}+N}s_1^{\lambda
_1}\cdots s_{q-1}^{\lambda _{q-1}}s_q^Nt^{\lambda _1+\cdots +\lambda
_{q-1}+N}}{\lambda _{1}!\cdots\lambda
_{q-1}!}T_0[T_{1}]^{[\lambda _{1}]}$$
$$\cdots[T_{q-1}]^{[\lambda _{q-1}]}\{ [T_q]^{[N]}(s_qt)\}
e^{-(s_{q+1}-s_q)tD^2}\cdots
T_{2k}e^{-(1-s_{2k})tF}.\eqno(4.16)$$
Since $T_0[T_{1}]^{[\lambda _{1}]}
\cdots[T_{q-1}]^{[\lambda _{q-1}]}$ is a $\lambda _{1}+\cdots+\lambda _{q-1}$ order differential operator and
 by Lemma 4.2 and Lemma 4.3, we get (see pp. 61-62 in [Fe])
$$||\psi_t\int_{\Delta_{2k}}t^k
\sum_{1\leq q\leq {2k}}\sum^{N-1}_{\lambda _{1},\cdots ,\lambda _{q-1}=0}
\frac{(-1)^{\lambda _1+\cdots +\lambda _{q-1}+N}s_1^{\lambda
_1}\cdots s_{q-1}^{\lambda _{q-1}}s_q^Nt^{\lambda _1+\cdots +\lambda
_{q-1}+N}}{\lambda _{1}!\cdots\lambda
_{q-1}!}T_0[T_{1}]^{[\lambda _{1}]}$$
$$\cdots[T_{q-1}]^{[\lambda _{q-1}]}\{ [T_q]^{[N]}(s_qt)\}
e^{-(s_{q+1}-s_q)tD^2}\cdots
T_{2k}e^{-(1-s_{2k})tF}dv||\sim O(t^{\frac{2k+N+\lambda_1+\cdots+\lambda_{2q}-{\rm dim}M}{2}}).\eqno(4.17)$$
So we get\\

\noindent{\bf Theorem 4.5}~~ (1)~{\it if~ $2k\leq {\rm dim}M$, then}
$${\rm Ch}_{2k}({\mathcal{B}}_t)(\phi)(h_0\otimes f_0,\cdots,h_{2k}\otimes f_{2k})~~~~~~~~~~~~~~~~~~~~~~~~~~~~~~~~~~~~~~~~~~~$$
$$=\psi_t
\sum^{{\rm dim}M-2k}_{\lambda_1,\dots,\lambda_{2k}=0}
\frac{
(-1)^{\lambda_1+\cdots+\lambda_{2k}}}{\lambda_1!\cdots \lambda _{2k}!}Ct^{|\lambda|+{k}}
{\rm
Str}[\widetilde{\phi}T_0[T_1]^{[\lambda_1]}\cdots [T_{2k}]^{[\lambda_{2k}]}e^{-tF}]+O(\sqrt {t}),\eqno(4.18)$$
{\it with the constant}\\
$$C=\frac{1}{\lambda_1+1}\frac{1}{\lambda_1+\lambda_2+2}\cdots \frac{1}{\lambda_1+\cdots +\lambda_{2k}+2k}.\eqno(4.19)$$
(2) {\it if~ $2k>{\rm dim}M$, then}
$${\rm Ch}_{2k}({\mathcal{B}}_t)(\phi)(h_0\otimes f_0,\cdots,h_{2k}\otimes f_{2k})=0.\eqno(4.20)$$\\

Write $P_{k,\lambda}:= T_0[T_1]^{[\lambda_1]}\cdots [T_{2k}]^{[\lambda_{2k}]}.$ We will compute
$${\rm lim}_{t\rightarrow 0^+}
\psi_tt^{|\lambda|+{k}}{\rm
Str}[\widetilde{\phi}P_{k,\lambda}e^{-tF}].$$
By Lemma 2.18 and Proposition 2.21, we have
$$\sigma(P_{k,\lambda})=T_0[F_{(2)},\sigma(T_1)]^{[\lambda_1]}\cdots [F_{(2)},\sigma(T_{2k})]^{[\lambda_{2k}]}+O_G(2|\lambda|+2k-1).\eqno(4.21)$$
Direct computations show that
$$[F_{(2)},\sigma(T_j)]=O_G(2),~~O_G([F_{(2)},\sigma(T_j)]^{[\lambda_{j}]})<2\lambda_j+1.\eqno(4.22)$$
Then when $(\lambda_1,\cdots,\lambda_{2k})\neq (0,\cdots, 0)$,
then $$O_G(P_{k,\lambda})=O_G(2|\lambda|+2k-1);~~O_G(P_{k,\lambda}(F+\partial_t)^{-1})=O_G(2|\lambda|+2k-3).\eqno(4.23)$$
By Lemma 3.5 and Lemma 3.6 (1), we have
$$\sigma[\psi_t I_{P_{k,\lambda}(F+\partial_t)^{-1}}(0,t)]^{(a,0)}=\omega^{\rm odd}O(t^{-|\lambda|-k+\frac{1}{2}})+O(t^{-|\lambda|-k+1}),
\eqno(4.24)$$
where $\omega^{\rm odd}\in \Omega^{\rm odd}(B).$
So in this case, $${\rm lim}_{t\rightarrow 0^+}
\psi_tt^{|\lambda|+{k}}{\rm
Str}[\widetilde{\phi}P_{k,\lambda}e^{-tF}]=0.\eqno(4.25)$$
When $(\lambda_1,\cdots,\lambda_{2k})=(0,\cdots, 0).$ Then $O_G(T_0T_1\cdots T_{2k})=2k$ and
$O_G(P_{k,\lambda}(F+\partial_t)^{-1})=O_G(2k-2)$. The model operator of $P_{k,\lambda}(F+\partial_t)^{-1}$ is $T_0\sigma(T_1)\cdots\sigma(T_{2k})(F_{(2)}+\partial_t)^{-1}$. By Lemma 3.6(2), we get
$${\rm lim}_{t\rightarrow 0^+}t^k\sigma[\psi_tI_{P_{k,\lambda}(F+\partial_t)^{-1}}(0,t)]^{(a,0)}
=T_0\sigma(T_1)\cdots \sigma(T_{2k})I_{(F_{(2)}+\partial_t)^{-1}}(0,1)^{(a,0)}.\eqno(4.26)$$
By (4.25) and (4.26) and the lemma 9.13 in [PW], we get\\

\noindent{\bf Theorem 4.6} {\it The following equality holds}
$${\rm lim}_{t\rightarrow 0^+}{\rm Ch}_{2k}({\mathcal{B}}_t)(\phi)(h_0\otimes f_0,\cdots,h_{2k}\otimes f_{2k})$$
$$
=\frac{1}{(2k)!}(-i)^{\frac{n}{2}}\sum_a
(2\pi)^{-\frac{a}{2}}\int_{M_a^\phi/B}T_0\sigma(T_1)\cdots\sigma(T_{2k})
\widehat{A}(R^{G^\phi}) \nu_\phi(R^{N^\phi}).\eqno(4.27)$$\\

\section { The regularity of the equivariant eta form}

{\bf 5.1  The regularity of the equivariant eta invariant}\\

\indent In [BiF], Bismut and Freed gave a simple proof of the regularity of
    eta invariants. In [Po2], Ponge gave another proof of the regularity of
    eta invariants using the method in [Po1].
      Zhang proved the regularity of the equivariant eta
    invariant by the Clifford asymptotics in [Zh3]. In this section, we shall extend the approach of [Po2] to the equivariant setting and prove the regularity of the
equivariant eta invariant by this approach.
Then we define the equivariant eta form and prove its regularity.\\
  \indent We shall give some notations. Let $X$ be a compact
oriented odd dimensional Riemannian manifold without boundary with a
fixed spin structure and $S$ be the bundle of spinor on $X$. Denote
by $D$ the associated Dirac operator on $H=L^2(X;S)$, the Hilbert
space of $L^2$-sections of the bundle $X$. Suppose that $\phi$ acts
on $X$ by orientation-preserving isometries and $\phi$ has a lift
$\widetilde{\phi}:~\Gamma(S)\rightarrow\Gamma(S)$ (see [LYZ]), then
we have $\widetilde{\phi}$ commutes with the Dirac operator and
$\widetilde{\phi}$ is a bounded operator. Then the equivariant eta
invariant is defined by
$$\eta_\phi(D)=\frac{1}{\sqrt{\pi}}\int_0^{\infty}\frac{1}{t^{\frac{1}{2}}}{\rm
Tr}[\widetilde{\phi}De^{-tD^2}]dt.\eqno(5.1)$$\\
We have\\

\noindent {\bf Theorem 5.1} ([Zh3]) {\it As $t\rightarrow 0^+$}
$${\rm
Tr}[\widetilde{\phi}De^{-tD^2}]\sim
O(t^{\frac{1}{2}}).\eqno(5.2)$$\\

~In order to prove Theorem 5.1, we introduce an auxiliary Grassmann variable $z$ as in [BiF], i. e. $z^2=0$. The auxiliary Grassmann variable $z$ may be considered as one form on the base
 $S^1$ of the fibre bundle $X\times S^1$ and set
$O_G(z)=1$. By the Duhamel principle, we have
$${\rm exp}(-t(D^2-zD))={\rm exp}(-tD^2)+ztD{\rm exp}(-tD^2).\eqno(5.3)$$
 Set
$$h(x)=1+\frac{1}{2}z\sum\limits^n_{i=1}{x_ic(e_i)},\eqno(5.4)$$
where $(x_1,\cdots,x_n)$ is the normal coordinates under the
parallel frame $e_1,\cdots, e_n$ and we consider $h$ as $h\chi$
where $\chi$ is a cut function about $(x_1,\cdots,x_n)$. By [Zh3], we
have
$$hc(e_i)h^{-1}=c(e_i)+O_G(0);~~h(D^2-zD)h^{-1}=D^2+zu,\eqno(5.5)$$
where $O_G(u)\leq 0$,~ $u$ contains no $z$ and the equality
$$ztD{\rm{exp}}(-tD^2)(x,y)=h^{-1}(x){\rm{exp}}(-t(D^2+zu)(x,y))h(y)-{\rm{exp}}(-tD^2)(x,y).\eqno(5.6)$$
Let
$$(D^2+zu)^{-1}=D^{-2}-zD^{-2}uD^{-2}.\eqno(5.7)$$
We may consider $D^2+zu$ as the operator with $1$-form coefficients
on $X\times S^1$. Then the Greiner's approach of the heat kernel
asymptotics with form coefficients in Section 2.1 works for
$D^2+zu$. By
$$(D^2+zu+\partial_t)^{-1}=(D^2+\partial_t)^{-1}-z(D^2+\partial_t)^{-1}u(D^2+\partial_t)^{-1}.\eqno(5.8)$$
$$K_{(D^2+zu+\partial_t)^{-1}}(x,y,t)=K_{(D^2+\partial_t)^{-1}}(x,y,t)-zK_{(D^2+\partial_t)^{-1}u(D^2+\partial_t)^{-1}}(x,y,t).\eqno(5.9)$$
Then
$$\widetilde{\phi}(x)h^{-1}(x)K_{(D^2+zu+\partial_t)^{-1}}(x,\phi(x),t)h(\phi(x))=
\widetilde{\phi}(x)h^{-1}(x)K_{(D^2+\partial_t)^{-1}}(x,\phi(x),t)h(\phi(x))$$
$$-z\widetilde{\phi}(x)K_{(D^2+\partial_t)^{-1}u(D^2+\partial_t)^{-1}}(x,\phi(x),t).\eqno(5.10)$$
Recall since the dimension $n$ is odd, we have
$${\rm Tr}[c(e^{i_1})\cdots c(e^{i_k})]=\left\{\begin{array}{cc}
2^{[n/2]} & ~~~~~~~~~{\rm if}~~k=0\\
0  &~~~~~~~~~~~~~~{\rm if }~0<k<n\\
(-i)^{[n/2]+1}2^{[n/2]} & ~~~~~~~~~{\rm if }~~~k=n
\end{array}\right.~~~~~~~\eqno(5.11)$$
Similar to Lemma 3.5, we have:\\

\noindent{\bf Lemma 5.2} {\it Let $A\in {\rm Cl}^{{\rm
odd}}(n)\subset{\rm End}(S_n)$, then}
$${\rm
Tr}[\widetilde{\phi}A]=(-i)^{\frac{n+1}{2}}2^{\frac{n-1}{2}}2^{-\frac{b}{2}}{\rm
det}^{\frac{1}{2}}(1-\phi^N)|\sigma(A)|^{(a,0)}$$
$$+(-i)^{\frac{n+1}{2}}2^{\frac{n-1}{2}}\sum_{0\leq b'<b}
|\sigma(\widetilde{\phi})^{(0,b')}\sigma(A)^{(a,b-b')}|^{(n)}.\eqno(5.12)$$\\

 We know
that Lemma 9.12 in [PW] also holds for the odd dimensional
manifold $X$ and when $a$ is odd. We know that
$(D^2+\partial_t)^{-1}u(D^2+\partial_t)^{-1}$ has the Getzler order $-4$
and odd Clifford elements. By Lemma 8.1 and Lemma 9.12 in [PW] and
Lemma 5.2 for $(D^2+\partial_t)^{-1}u(D^2+\partial_t)^{-1}$ and
$j\geq a$, so
$$\int_X{\rm
Tr}[\widetilde{\phi}K_{(D^2+\partial_t)^{-1}u(D^2+\partial_t)^{-1}}(x,\phi(x),t)]dx=O(t^{\frac{3}{2}}).\eqno(5.13)$$
By equalities (5.10) and (5.13), in order to prove Theorem 5.1, we need to prove
$$\widetilde{\phi}(x)h^{-1}(x)K_{(D^2+\partial_t)^{-1}}(x,\phi(x),t)h(\phi(x))
=\widetilde{\phi}(x)K_{(D^2+\partial_t)^{-1}}(x,\phi(x),t)+O(t^{\frac{3}{2}}).\eqno(5.14)$$
This comes from the equality
$$h^{-1}(x',v)c(e_i)h(x',\phi'(x)v)=c(e_i)+\frac{1}{2}zc(e_i)\sum_{j=a+1}^n[(\phi'(x)-1)v]_jc(e_j)-zx_i,\eqno(5.15)$$
and $\widetilde{\phi}$ only contains the Clifford elements
$c(e_{a+1}),\cdots,c(e_n)$.\\

\noindent{\bf 5.2~The regularity of the equivariant eta form}\\

 In this section, the fundamental setup is the same as Section
 2.1. But we assume the dimension $n$ of the fibre to be odd.  Let $\phi$ be an isometry which acts fiberwise on $M$. We will consider
that $\phi$ acts as identity on $B$. Let ${\rm Tr}^{\rm even}$
be the trace on the coefficients of even forms on $B$. We assume that ${\rm ker}D$ is a complex vector bundle. Then
the equivariant eta form is defined by
$$\widehat{\eta}(\phi):=\int_0^{\infty}{\rm Tr}^{\rm
even}[\widetilde{\phi}\frac{d{\mathcal{B}}_t}{dt}e^{-{\mathcal{B}}_t^2}]dt.\eqno(5.16)$$
Let  $T$ be the torsion tensor of $\nabla^{\oplus}$ and
$c(T)=\sum_{1\leq \alpha<\beta \leq m}dy_\alpha dy_\beta
c(T(\frac{\partial}{\partial y_\alpha},\frac{\partial}{\partial
y_\beta})).$ Then
$$\widehat{\eta}(\phi)=\int_0^{\infty}\frac{1}{2\sqrt{t}}{\rm Tr}^{\rm
even}[\psi_t\widetilde{\phi}(D+\frac{c(T)}{4})e^{-tF}]dt.\eqno(5.17)$$\\

\noindent {\bf Lemma 5.3} {\it We assume that ${\rm ker}D$ is a complex vector bundle. As $t\rightarrow +\infty$, we have}
$$\frac{1}{2\sqrt{t}}{\rm Tr}^{\rm
even}[\psi_t\widetilde{\phi}(D+\frac{c(T)}{4})e^{-tF}]\sim
O(t^{-\frac{3}{2}}).\eqno(5.18)$$\\

\noindent{\bf Proof.} Since $\widetilde{\phi}$ is a bounded operator, so the proof of Lemma 5.3 is the same as the proof of Theorem 9.7 in [BGV].  $\Box$\\

 \indent Let $e_1(x),\cdots e_n(x)$ denote the orthonormal frame of $TZ$.
   If $A(Y)$ is any $0$ order operator depending linearly on
   $Y\in \Gamma (M,TM)$, we define the operator $(\nabla_{e_i}+A(e_i))^2$ as
   follows
   $$(\nabla_{e_i}+A(e_i))^2=\sum_1^n(\nabla_{e_i(x)}+A(e_i(x)))^2-\nabla_{\sum_j\nabla_{e_j}e_j}
   -A(\sum_j\nabla_{e_j}e_j).\eqno(5.21)$$
 Then by [Bi2], we have
$$F=-(\nabla^{TZ}_{e_i}+\frac{1}{2}<S(e_i)e_j,
f_\alpha>e_jdy_\alpha+\frac{1}{4}<S(e_i)f_\alpha,f_\beta>dy_\alpha
dy_\beta)^2+\frac{r}{4}.\eqno(5.22)$$
 Set
$$h(x)=1+\frac{1}{2}dt\sum_{i=1}^nx_ie_i.\eqno(5.23)$$
Then we have
$$ he_ih^{-1}=e_i+dtL_1;~~
h(\frac{1}{2}e_idt)h^{-1}=\frac{1}{2}e_idt,\eqno(5.24)$$
$$~h\nabla^G_{e_i}h^{-1}=\nabla^G_{e_i}-\frac{1}{2}dte_i+dtL_2,~ \eqno(5.25)$$
$$h(\frac{1}{2}<S(e_i)e_j,
f_\alpha>e_jdy_\alpha )h^{-1}=\frac{1}{2}<S(e_i)e_j,
f_\alpha>e_jdy_\alpha +dtL_3,\eqno(5.26)$$
$$h(\frac{1}{4}<S(e_i)f_\alpha,f_\beta>dy_\alpha
dy_\beta)h^{-1}=\frac{1}{4}<S(e_i)f_\alpha,f_\beta>dy_\alpha
dy_\beta,\eqno(5.27)$$
where $L_1,L_2,L_3\in O_G(-1).$ By Proposition 2.10 in [BeGS], we have
$$F+dt(D+\frac{c(T)}{4})=-(\nabla^{TZ}_{e_i}+\frac{1}{2}\left<S(e_i)e_j,
f_\alpha\right>e_jdy_\alpha$$
$$+\frac{1}{4}\left<S(e_i)f_\alpha,f_\beta\right>dy_\alpha
dy_\beta-\frac{1}{2}e_idt)^2+\frac{r}{4},\eqno(5.28)$$ By
(5.24)-(5.28), we have
$$h[F+dt(D+\frac{c(T)}{4})]h^{-1}=F+dtu,\eqno(5.29)$$
where $O_G(u)\leq 0$.  By the Duhamel principle, we can get
$$tdt{\rm
Tr}^{\rm even}[\widetilde{\phi}(D+\frac{c(T)}{4}){\rm exp}(-tF)]
={\rm Tr}^{\rm even}[\widetilde{\phi}{\rm exp}(-tF)]-{\rm Tr}^{\rm
even}[\widetilde{\phi}h^{-1}{\rm exp}(-t(F+dtu))h].\eqno(5.30)$$
Then
$$t\frac{dt}{\sqrt{t}}{\rm
Tr}^{\rm even}[\widetilde{\phi}(D+\frac{c(T)}{4t})\psi_t{\rm
exp}(-tF)] ={\rm Tr}^{\rm even}[\widetilde{\phi}\psi_t{\rm
exp}(-tF)]$$ $$-{\rm Tr}^{\rm even}[\widetilde{\phi}h^{-1}\psi_t{\rm
exp}(-t(F+dtu))h].\eqno(5.31)$$ We know that Lemma 3.6 is still
correct. When we take ${\rm Tr}^{{\rm even}}$, the terms having
coefficients in $\Omega^{\rm odd}$ vanish. Similar to Section
5.1, we can get\\

\noindent {\bf Theorem 5.4} {\it As $t\rightarrow 0^+$, we have}
$${\rm Tr}^{\rm
even}[\psi_t\widetilde{\phi}(D+\frac{c(T)}{4})e^{-tF}]\sim
O(t^{\frac{1}{2}}).\eqno(5.32)$$\\

 By Lemma 5.3 and Theorem 5.4, the equivariant eta form is well
 defined. When $n$ is even, we define
$$\widehat{\eta}(\phi):=\int_0^{\infty}{\rm Str}[\widetilde{\phi}\frac{d{\mathcal{B}}_t}{dt}e^{-{\mathcal{B}}_t^2}]dt.\eqno(5.33)$$
Similarly, we may prove the regularity in this case.\\

 \noindent {\bf Acknowledgement.} This project
was supported by NSFC No.11271062 and NCET-13-0721. The author would like to thank Profs. Weiping Zhang
and Huitao Feng for introducing index theory to him. The author would
like to thank Prof. Weiping Zhang and Dr. Hang Wang for helpful discussions and comments.  The author would
like to thank referees for careful reading and helpful comments.\\

\noindent{\large \bf References}\\

\noindent[APS]M. F. Atiyah; V. K. Patodi; I. M. Singer, Spectral
asymmetry and Riemannian geometry, I. Math. Proc. Cambridge Philos.
Soc. 77(1975), 43-69.

\noindent[AS]M. F. Atiyah; I. M. Singer, The index of elliptic
operators IV, Annals of Math. 93(1971), 119-138.

\noindent[Az1]F. Azmi, The equivariant Dirac cyclic cocycle, Rocky Mountain J. Math. 30(2000), 1171-1206.

\noindent[Az2]F. Azmi, The Chern-Connes character formula for families of Dirac operators,
Publ. Math. Debrecen 61(2002), 1-27.

\noindent[Az3]F. Azmi, Equivariant bivariant cyclic theory and equivariant Chern-Connes character, Rocky Mountain J. Math. 34(2004), no. 2, 391-412.

\noindent[BeGS]R. Beals; P. Greiner; N. Stanton, The heat equation on a CR manifold, J. Differential Geom. 20(1984),
343-387.

\noindent[BeC]M. Benameur; A. Carey, Higher spectral flow and an entire bivariant JLO cocycle, J. K-Theory 11(2013), no. 1, 183-232.

\noindent[BGV]N. Berline; E. Getzler; M. Vergne, {\it Heat kernels
and Dirac operators}, Springer-Verlag, Berlin, 1992.

\noindent[BV1]N. Berline; M. Vergne, A proof of Bismut local index theorem for a family of Dirac operators, Topology 26(1987), no. 4, 435-463.

\noindent[BV2]N. Berline; M. Vergne, A computation of the
equivariant index of the Dirac operators, Bull. Soc. Math. Prance
113(1985) 305-345.

 \noindent[Bi1]J.-M. Bismut, The Atiyah-Singer theorems: A
probabilistic approach, J. Func. Anal. 57(1984) 56-99.

\noindent[Bi2]J. M. Bismut, The Atiyah-Singer index theorem for
families of Dirac operators: Two heat equation proofs, Invent. Math.
83(1986) 91-151.

\noindent[Bi3]J. M. Bismut, Equivariant immersions and Quillen metrics, J. Diff. Geom. 41
(1995). 53-159

\noindent[BC1]J. M. Bismut; J. Cheeger, Families index for
manifolds with boundary, superconnections, and cones. I. Families of
manifolds with boundary and Dirac operators. J. Funct. Anal. 89
(1990), no.2, 313-363.

\noindent[BC2]J. M. Bismut; J. Cheeger, Families index for
manifolds with boundary, superconnections and cones. II. The Chern
character, J. Funct. Anal. 90(1990), no. 2, 306-354.

\noindent[BiF]J. M. Bismut; D. Freed, The analysis of elliptic
families. II. Dirac operators, eta invariants, and the holonomy
theorem, Comm. Math. Phys. 107(1986), no. 1, 103-163.

\noindent[BiGS]J. M .Bismut; H. Gillet; C. Soul$\acute{e}$, Analytic torsion
and holomorphic determinant bundles. I. Bott-Chern forms and
analytic torsion, Comm. Math. Phys. 115(1988), no. 1, 49-78.

\noindent[BlF]J. Block; J. Fox, Asymptotic pseudodifferential operators and index theory, Contemp. Math., 105(1990), 1-32.

\noindent[CH]S. Chern; X. Hu, Equivariant Chern character
for the invariant Dirac operators, Michigan Math. J. 44(1997),
451-473.

\noindent[Co]A. Connes, Entire cyclic cohomology of Banach algebras and the character of
¦È-summable Fredholm module, K-Theory 1(1988), 519-548.

\noindent[CM1]A, Connes; H. Moscovici, Cyclic cohomology, the Novikov conjecture and hyperbolic groups,
Topology 29(1990), no. 3, 345-388.

\noindent [CM2]A. Connes; H. Moscovici, Transgression and
Chern character of finite dimensional K-cycles,
 Commun. Math. Phys. 155(1993), 103-122.

\noindent[CFKS]H. L. Cycon; R. G. Froese; W. Kirsch; B. Simon, {\it Schrodinger Operators
with Application to Quantum Mechanics and Global Geometry}, Springer-Verlag, Texts and Monographs in Physics, 1987.

\noindent[Do1]H. Donnelly, Local index theorem for families, Michigan Math. J. 35(1988), no. 1, 11-20.

\noindent[Do2]H. Donnelly, Eta invariants for $G$-spaces, Indiana
Univ. Math. J. 27 (1978), no. 6, 889-918.

\noindent[Fe]H. Feng, A note on the noncommutative Chern
character (in Chinese), Acta Math. Sinica 46 (2003), 57-64.

\noindent[Ge1]E. Getzler, PseudodifferentiaI operators on
supermanifolds and the Atiyah-Singer index theorem, Commun. Math.
Phys. 92(1983), 163-178.

\noindent[Ge2]E. Getzler, A short proof of the local Atiyah-Singer
index theorem, Topology 25(1986), 111-117.

\noindent[Gr]P. Greiner, An asymptotic expansion for the heat equation, Arch. Rational Mech. Anal. 41(1971), 163-218.

\noindent[JLO]A. Jaffe, A. Lesniewski; K. Osterwalder, Quantum K-theory: The Chern
character, Comm. Math. Phys. 118(1988), 1-14.

\noindent [LYZ]J. D. Lafferty, Y. L. Yu; W. P. Zhang, A
direct geometric proof of Lefschetz fixed point formulas, Trans.
AMS. 329 (1992), 571-583.

\noindent [LM]K. Liu; X. Ma, On family rigidity theorems. I. Duke Math. J. 102(2000), no. 3, 451-474.

\noindent[Lo]J. Lott, Superconnections and higher index theory, Geom. Func. Anal. 2
(1992), 421-454.

\noindent[MS]H. P. Mckean; I. M. Singer, Curvature and the
eigenvalue of the Laplacian, J. Diff. Geom., I(1967), 43-69.

 \noindent[Pa]V. K. Patodi, Curvatures and the eigenforms of the Laplace operator, J. Diff. Geom., 5(1971), 238-249.

\noindent[Po1]R. Ponge, A new short proof of the local index
formula and some of its applications, Comm. Math. Phys. 241(2003), 215-234.

\noindent[Po2]R. Ponge, A new proof of the local regularity of the eta invariant of a Dirac operator,
J. Geom. Phys. 56(2006), no. 9, 1654-1665.

\noindent[PW]R. Ponge; H. Wang, Noncommutative geometry, conformal geometry, and the local equivariant index theorem,
arXiv:1411.3703.

\noindent[Qu]D. Quillen, Superconnections and the Chern character,
Topology 24(1985) 89-95.

\noindent[Si]B. Simon, {\it Trace ideals and their applications}, London Math. Soc., Lecture Note 35, 1979.

\noindent[Ta]M. E. Taylor, {\it Pseudodifferential operators}, Princeton Mathematical Series, 34. Princeton University Press, Princeton, N.J., 1981.

 \noindent[Wu]F. Wu, A bivariant Chern-Connes character and the higher $\Gamma$-index theorem, K-Theory 11(1997), no. 1, 35-82.

 \noindent[Yu]Y. Yu, Local index theorem for Dirac operators,
Acta. Math. Sinica., 3(1987), No.2, 152-169.

 \noindent[Zh1]W. Zhang, Local Atiyah-Singer index
theorem for families of Dirac operators, Differential geometry and
topology, Lecture Notes in Math., 1369, Springer, Berlin, 1989, 351-366.

\noindent[Zh2]W. Zhang, The Chern-Connes character and the local family index theorem, preprint.

\noindent[Zh3]W. Zhang, A note on equivariant eta invariants,
Proc. AMS 108(1990), no. 4, 1121-1129.\\

 \indent{\it School of Mathematics and Statistics, Northeast Normal University, Changchun Jilin, 130024, China }\\
 \indent E-mail: {\it wangy581@nenu.edu.cn}\\

\end{document}